# SUFFICIENT CONDITIONS ON THE EXISTENCE OF SWITCHING OBSERVERS FOR NONLINEAR TIME-VARYING SYSTEMS

D. BOSKOS [*] and J. TSINIAS [**]

**Abstract.** We derive sufficient conditions for the solvability of the observer design problem for a wide class of nonlinear time-varying systems, including those having triangular structure. We establish that, under weaker assumptions than those imposed in the existing works in the literature, it is possible to construct a switching sequence of time-varying noncausal dynamics, exhibiting the state determination of our system.

**Matematics Subject Classifications.** 93B30, 93B51, 93C10

## I. INTRODUCTION

Several important contributions towards solvability of the Observer Design Problem (ODP) have been appeared in the literature; see for instance [1]-[26],[28]-[51]. In the present work we generalize the Lyapunov like approach of the works [46],[47],[48], and particularly of the recent paper [49], in order to derive sufficient conditions for the solvability of the ODP for a wide class of nonlinear systems under weaker assumptions than those imposed in the existing works on the same subject. The main idea of our work is to construct a switching sequence of time-varying noncausal dynamics exhibiting the state determination of our system. It should be mentioned that the idea of using switching observers has been adopted in several earlier works; see for instance [3],[17],[18],[36]. The results of the present work generalize those obtained in the previously mentioned contributions.

We consider time-varying systems of the form:

$$\dot{x} = f(t,x) \quad (1.1a)$$
$$y = h(t,x) \quad (1.1b)$$
$$(t,x) \in \mathbb{R}_{\geq 0} \times \mathbb{R}^n, \ y \in \mathbb{R}^k$$

where $y(\cdot)$ plays the role of output. In Section II we provide sufficient conditions for the solvability of the ODP for the general case (1.1) with linear output dynamics by means of a *non-causal* observer. The corresponding results of this section (Propositions 2.2 and 2.3) constitute extensions of Proposition 2.1 in [49]. In Section III (Proposition 3.1 and Corollary 3.1) we use the results of Section II for the derivation of sufficient conditions for the solvability of the ODP for a class of composite systems of the following form:

$$\dot{x}_1 = f_1(t,x_1) + G(t,x_1)x_2$$
$$\dot{x}_2 = f_2(t,x_1,x_2), \ (x_1,x_2) \in \mathbb{R}^{n_1} \times \mathbb{R}^{n_2} \quad (1.2a)$$
$$y = x_1 \quad (1.2b)$$

and in Proposition 4.1 of Section IV we exploit the results obtained in Sections II and III to establish sufficient conditions for the solvability of the ODP for triangular systems:

___________

Department of Mathematics, National Technical University of Athens, Zografou Campus 15780, Athens, Greece
e-mail adresses: dmposkos@central.ntua.gr [*], jtsin@central.ntua.gr [**]

$$\dot{x}_i = f_i(t, x_1, x_2, ..., x_i) + a_i(t, x_1)x_{i+1}, \; i = 1, 2, ..., n-1$$
$$\dot{x}_n = f_n(t, x_1, x_2, ..., x_n), \; (x_1, x_2, ..., x_n) \in \mathbb{R}^n \quad (1.3a)$$
$$y = x_1 \quad (1.3b)$$

The corresponding result generalizes Proposition 3.1 in [49]. Finally, in Section V we derive sufficient conditions for the solvability of the ODP for the general case (1.1), as well for (1.2) and (1.3), by means of a *causal* observer (Proposition 5.3, Corollary 5.2 and Proposition 5.5).

*Notations and Definitions*

We adopt the following notations. For a given vector $x \in \mathbb{R}^n$, $x'$ denotes its transpose and $|x|$ its Euclidean norm. We use the notation $|A| := \max\{|Ax| : x \in \mathbb{R}^n; |x| = 1\}$ for the induced norm of a matrix $A \in \mathbb{R}^{m \times n}$. By $N$ we denote the class of all increasing $C^0$ functions $\phi : \mathbb{R}_{\geq 0} \to \mathbb{R}_{\geq 0}$. For given $R > 0$, we denote by $B_R$ the closed ball of radius $R$, centered at $0 \in \mathbb{R}^n$. Let $t_0 \geq 0$ and consider an integer $\ell \in \mathbb{N}$ and a set-valued map $[t_0, \infty) \ni t \to Q(t) \subset \mathbb{R}^\ell$. We say that $Q(\cdot)$ satisfies the *Compactness Property* (**CP**), if for every sequence $(t_\nu)_{\nu \in \mathbb{N}} \subset [t_0, \infty)$ with $t_\nu \to t \in [t_0, \infty)$ and $q_\nu \in Q(t_\nu)$ there exist a subsequence $(q_{\nu_k})_{k \in \mathbb{N}}$ and a vector $q \in Q(t)$ such that $q_{\nu_k} \to q$.

**Definition 1.1.** *Let $k, \ell, n$ be positive integers, let $M$ and $S$ be nonempty subsets of $\mathbb{R}^n$ and $\mathbb{R}^\ell$, respectively, and let $\Omega(\mathbb{R}_{\geq 0}, M)$ be a nonempty set of functions $y = y_{t_0, x_0} : [t_0, \infty) \to \mathbb{R}^k$ parameterized by $(t_0, x_0) \in \mathbb{R}_{\geq 0} \times M$. For given $0 \leq \tau \leq \infty$ we say that the map*

$$([t_0, \infty) \times \mathbb{R}^n) \times (\mathbb{R}_{\geq 0} \times \Omega(\mathbb{R}_{\geq 0}, M)) \ni (t, x; t_0, y) \to a_{t_0, y}(t, x) \in S$$

*is $\tau$-noncausal with respect to $\Omega(\mathbb{R}_{\geq 0}, M)$, if for any $x \in \mathbb{R}^n$, $t \geq t_0 \geq 0$ and $y \in \Omega(\mathbb{R}_{\geq 0}, M) \cap C^0([t_0, \infty); \mathbb{R}^k)$ the value $a(t) := a_{t_0, y}(t, x)$ depends only on the values $\{y(s) : t \leq s \leq t + \tau\}$ of $y(\cdot)$. We say that $a_{t_0, y}(\cdot, \cdot)$ is causal with respect to $\Omega(\mathbb{R}_{\geq 0}, M)$, if it is 0-noncausal with respect to $\Omega(\mathbb{R}_{\geq 0}, M)$, namely, for any $x \in \mathbb{R}^n$, $t \geq t_0 \geq 0$ and $y \in \Omega(\mathbb{R}_{\geq 0}, M) \cap C^0([t_0, \infty); \mathbb{R}^k)$ the value $a(t) := a_{t_0, y}(t, x)$ is independent of the future values $\{y(s) : s > t\}$ of $y(\cdot)$.*

**Remark 1.1.** Obviously, if $a_{t_0, y}(\cdot, \cdot)$ is $\tau_0$-noncausal with respect to $\Omega(\mathbb{R}_{\geq 0}, M)$, then it is $\tau$-noncausal with respect to $\Omega(\mathbb{R}_{\geq 0}, M)$ for every $\tau \geq \tau_0$. Also, if $\varnothing \neq D \subset M$ and $y \in \Omega(\mathbb{R}_{\geq 0}, D)$ then, if the map $a_{t_0, y}(\cdot, \cdot)$ is $\tau_0$-noncausal with respect to $\Omega(\mathbb{R}_{\geq 0}, M)$, it is $\tau_0$-noncausal with respect to $\Omega(\mathbb{R}_{\geq 0}, D)$ as well.

**Definition 1.2.** *Let $\varnothing \neq M \subset \mathbb{R}^n$. We say that (1.1a) is M-forward complete, if there exists a function $\beta \in NNN$ such that the solution $x(t) := x(t, t_0, x_0)$ of (1.1a) initiated from $x_0$ at time $t = t_0$ satisfies:*

$$|x(t)| \leq \beta(t, t_0, |x_0|), \; \forall t \geq t_0 \geq 0, \; x_0 \in M \quad (1.4)$$

**Remark 1.2.** For the case $M = \mathbb{R}^n$ it is known (see, Lemma 2.3 in [27]) that existence of the solution $x(\cdot, t_0, x_0)$ of (1.1a) for all $t \geq t_0 \geq 0$ and $x_0 \in \mathbb{R}^n$ is equivalent to the existence of a function $\beta \in NNN$ satisfying (1.4).

In the sequel, for any pair of nonempty sets $I$, $M$ of $\mathbb{R}_{\geq 0}$ and $\mathbb{R}^n$, respectively, for which (1.1a) is *M*-forward complete, we adopt the notation $Y(I, M)$ to denote the set of all outputs of (1.1) with initial $(t_0, x_0) \in I \times M$, namely:

$$Y(I,M) := \{y : [t_0, \infty) \to \mathbb{R}^k : y(t) = h(t, x(t, t_0, x_0))\, ; \, (t_0, x_0) \in I \times M\} \tag{1.5}$$

**Definition 1.3.** *Let $\emptyset \neq M \subset \mathbb{R}^n$ and assume that (1.1a) is M-forward complete. We say that the Almost Causal Observer Design Problem (AC-ODP) is solvable for (1.1) with respect to $Y(\mathbb{R}_{\geq 0}, M)$, if for every $\bar{t}_0 \geq t_0 \geq 0$, $\tau > 0$ and $y \in Y(\mathbb{R}_{\geq 0}, M) \cap C^0([t_0, \infty); \mathbb{R}^k)$ there exist a continuous map:*

$$g := g_{\bar{t}_0, y}(t, z, w) : [\bar{t}_0, \infty) \times \mathbb{R}^n \times \mathbb{R}^k \to \mathbb{R}^n \tag{1.6}$$

*being $\tau$-noncausal with respect to $Y(\mathbb{R}_{\geq 0}, M)$, and a nonempty set $\bar{M} \subset \mathbb{R}^n$ such that for every $z_0 \in \bar{M}$ the corresponding trajectory $z(\cdot) := z(\cdot, \bar{t}_0, z_0)$ ; $z(\bar{t}_0) = z_0$ of the observer:*

$$\dot{z} = g(t, z, y) \tag{1.7}$$

*exists for all $t \geq \bar{t}_0$ and the error $e(t) := x(t) - z(t)$ between the trajectory $x(\cdot) := x(\cdot, t_0, x_0)$, $x_0 \in M$ of (1.1a) and the trajectory $z(\cdot) := z(\cdot, \bar{t}_0, z_0)$ of (1.7) satisfies:*

$$\lim_{t \to \infty} e(t) = 0 \tag{1.8}$$

*We say that the Strong Observer Design Problem (S-ODP) is solvable for (1.1) with respect to $Y(\mathbb{R}_{\geq 0}, M)$, if, in addition to (1.8), the map $g(\cdot)$ in (1.6) is causal with respect to $Y(\mathbb{R}_{\geq 0}, M)$. We say that the Almost Causal Switching Observer Design Problem (AC-SODP) is solvable for (1.1) with respect to $Y(\mathbb{R}_{\geq 0}, M)$, if for every $t_0 \geq 0$, $\tau > 0$ and $y \in Y(\mathbb{R}_{\geq 0}, M) \cap C^0([t_0, \infty); \mathbb{R}^k)$ there exist:*

- *a strictly increasing sequence of times $(t_m)_{m \in \mathbb{N}}$ with*

$$t_1 = t_0 \text{ and } \lim_{m \to \infty} t_m = \infty \tag{1.9}$$

- *a sequence of continuous mappings*

$$g_m := g_{m, t_m, y}(t, z, w) : [t_{m-1}, t_{m+1}] \times \mathbb{R}^n \times \mathbb{R}^k \to \mathbb{R}^n,\ m \in \mathbb{N} \tag{1.10}$$

*being $\tau$-noncausal with respect to $Y(\mathbb{R}_{\geq 0}, M)$, and a nonempty set $\bar{M} \subset \mathbb{R}^n$ such that any solution $z_m(\cdot)$ of the system*

$$\dot{z}_m = g_m(t, z_m, y),\ t \in [t_{m-1}, t_{m+1}],\ m \in \mathbb{N} \tag{1.11}$$

*with initial $z(t_{m-1}) \in \bar{M}$ is defined for each $t \in [t_{m-1}, t_{m+1}]$ and in such a way that, if we consider the piecewise continuous map $Z : [t_0, \infty) \to \mathbb{R}^n$ defined as:*

$$Z(t) := z_m(t),\ t \in [t_m, t_{m+1}),\ m \in \mathbb{N} \tag{1.12}$$

*where for each $m \in \mathbb{N}$ the map $z_m(\cdot)$ denotes the solution of (1.11), then the error $e(t) := x(t) - Z(t)$ between the trajectory $x(\cdot) := x(\cdot, t_0, x_0)$, $x_0 \in M$ of (1.1a) and $Z(\cdot)$ satisfies (1.8). Finally, we say that the Strong Switching Observer Design Problem (S-SODP) is solvable for (1.1) with respect to $Y(\mathbb{R}_{\geq 0}, M)$, if, in addition to (1.8), for each $m \in \mathbb{N}$ the map $g_m(\cdot)$ in (1.10) is causal with respect to $Y(\mathbb{R}_{\geq 0}, M)$.*

**Remark 1.3.** As mentioned in [49], solvability of AC-(S)ODP is equivalent to the solvability of the strong ODP by means of a *time-delay* system with *causal* dynamics.

## II. THE GENERAL CASE

In this section we specialize our analysis to the case of systems with linear output:

$$\dot{x} = f(t,x) := F(t, x, H(t)x), \quad (t,x) \in \mathbb{R}_{\geq 0} \times \mathbb{R}^n \tag{2.1a}$$

$$y = h(t,x) := H(t)x, \quad y \in \mathbb{R}^k \tag{2.1b}$$

where $H : \mathbb{R}_{\geq 0} \to \mathbb{R}^{k \times n}$ is $C^0$ and $F : \mathbb{R}_{\geq 0} \times \mathbb{R}^n \times \mathbb{R}^k \to \mathbb{R}^n$ is $C^0$ and locally Lipschitz on $(x,y) \in \mathbb{R}^n \times \mathbb{R}^k$. We also assume that there exists a nonempty subset $M$ of $\mathbb{R}^n$ with $0 \in \mathrm{cl}\, M$ such that system (2.1a) is $M$-forward complete, namely, the solution $x(t) := x(t, t_0, x_0)$ of (2.1a) satisfies (1.4) for certain $\beta \in NNN$. In addition to (1.4), we make the following hypotheses. We assume that there exist an integer $\ell \in \mathbb{N}$, a map

$$A \in C^0(\mathbb{R}_{\geq 0} \times \mathbb{R}^\ell \times \mathbb{R}^k ; \mathbb{R}^{n \times n}) \tag{2.2}$$

and constants $L > 1$ and $R > 0$ such that the following properties hold:

**A1** *For every $t_0 \geq 0$ and $\xi > 0$ there exists a set-valued map*

$$[t_0, \infty) \ni t \to Q_R(t) := Q_{R, t_0, \xi}(t) \subset \mathbb{R}^\ell \tag{2.3}$$

*satisfying the CP (see notations) in such a way that for every $t \geq t_0$ and for every*

$$y \in Y_R(t) := \{ y \in \mathbb{R}^k : y = H(t)x, \, |x| \leq \beta(t, t_0, R) \}, \tag{2.4}$$
$$x, z \in \mathbb{R}^n \text{ with } |x| \leq \beta(t, t_0, R) \text{ and } |x - z| \leq \xi$$

*the following holds:*

$$\Delta F(t, x, z; y) := F(t, x, y) - F(t, z, y)$$
$$= A(t, q, y)(x - z) \text{ for certain } q \in Q_R(t) \tag{2.5}$$

**A2** *There exists a constant $\varepsilon_R > 0$, (being independent of the constant $L$) such that for every $t_0 \geq 0$ and $\xi > 0$ there exists a set-valued map $Q_R := Q_{R, t_0, \xi}$ as in (2.3) satisfying the CP in such a way that, for every $\bar{t}_0 \geq t_0$, $\tau_0 > 0$ and $y \in Y(\mathbb{R}_{\geq 0}, M) \cap C^0([t_0, \infty); \mathbb{R}^k)$, a time-varying symmetric matrix*

$$P_R := P_{R, t_0, \bar{t}_0, \tau_0, \xi, y} \in C^1([\bar{t}_0, \infty); \mathbb{R}^{n \times n})$$

*and a function*

$$d_R := d_{R, t_0, \bar{t}_0, \tau_0, \xi, y} \in C^0([\bar{t}_0, \infty); \mathbb{R})$$

*can be found, both being $\tau_0$-noncausal with respect to $Y(\mathbb{R}_{\geq 0}, M)$, satisfying:*

$$P_R(t) \geq I_{n \times n}, \quad \forall t \geq \bar{t}_0 \,;\, |P_R(\bar{t}_0)| \leq L; \tag{2.6a}$$

$$\int_{\bar{t}_0}^t d_R(s) ds > -\varepsilon_R, \quad \forall t \geq \bar{t}_0 \,;\, \int_{\bar{t}_0}^\infty d_R(s) ds = \infty; \tag{2.6b}$$

*moreover, the following holds*

$$e' P_R(t) A(t, q, y(t)) e + \tfrac{1}{2} e' \dot{P}_R(t) e \leq -d_R(t) e' P_R(t) e, \tag{2.6c}$$
$$\forall t \geq \bar{t}_0, \, e \in \ker H(t), \, q \in Q_R(t)$$

*provided that*

$$y \in Y(\mathbb{R}_{\geq 0}, B_R \cap M) \cap C^0([t_0, \infty); \mathbb{R}^k) \tag{2.7}$$

We show that under A1 and A2 the AC-ODP is solvable for the case (2.1) with respect to $Y(\mathbb{R}_{\geq 0}, B_R \cap M)$. We first need to establish a technical preliminary result. Let $k, \ell, n, p \in \mathbb{N}$, consider a pair $(H, A)$ of continuous mappings:

$$H := H(t, y) \,;\; H : \mathbb{R}_{\geq 0} \times \mathbb{R}^k \to \mathbb{R}^{p \times n}; \tag{2.8a}$$

$$A := A(t, q, y) \,;\; A : \mathbb{R}_{\geq 0} \times \mathbb{R}^\ell \times \mathbb{R}^k \to \mathbb{R}^{n \times n}, \tag{2.8b}$$

let $U$ and $W$ be nonempty subsets of $\mathbb{R}^n$ with $U \cap W \neq \emptyset$ and let $\Omega(\mathbb{R}_{\geq 0}, U)$ be a (nonempty) set of functions $y := y_{t_0, x_0} : [t_0, \infty) \to \mathbb{R}^k$ parameterized by $(t_0, x_0) \in \mathbb{R}_{\geq 0} \times U$. We make the following hypothesis for the pair $(H, A)$ above:

**H1** *For every $t_0 \geq 0$ there exists a set-valued map $[t_0, \infty) \ni t \to Q(t) \subset \mathbb{R}^\ell$ satisfying the CP in such a way that for every $\bar{t}_0 \geq t_0$, $\tau_0 > 0$ and $y \in \Omega(\mathbb{R}_{\geq 0}, U) \cap C^0([t_0, \infty); \mathbb{R}^k)$ there exist a time-varying positive definite matrix $P \in C^1([\bar{t}_0, \infty); \mathbb{R}^{n \times n})$ and a function $d \in C^0([\bar{t}_0, \infty); \mathbb{R})$, both $\tau_0$-noncausal with respect to $\Omega(\mathbb{R}_{\geq 0}, U)$, such that*

$$\begin{aligned} &e'P(t)A(t,q,y(t))e + \tfrac{1}{2}e'\dot{P}(t)e \leq -d(t)e'P(t)e, \\ &\forall t \geq \bar{t}_0,\, e \in \ker H(t, y(t)),\, q \in Q(t) \end{aligned} \tag{2.9}$$

*provided that*

$$y \in Y(\mathbb{R}_{\geq 0}, U \cap W) \cap C^0([t_0, \infty); \mathbb{R}^k) \tag{2.10}$$

**Proposition 2.1.** *Consider the pair $(H, A)$ of the continuous time-varying mappings in (2.8a,b) and assume that H1 is fulfilled. Then for every $\bar{t}_0 \geq t_0 \geq 0$, $\tau > \tau_0 > 0$, $y \in \Omega(\mathbb{R}_{\geq 0}, U) \cap C^0([t_0, \infty); \mathbb{R}^k)$ and $\bar{d} \in C^0([\bar{t}_0, \infty); \mathbb{R})$ with*

$$\bar{d}(t) < d(t), \quad \forall t \geq \bar{t}_0 \tag{2.11}$$

*there exists a function $\phi = \phi_{t_0, \bar{t}_0, \tau, y} \in C^1([\bar{t}_0, \infty); \mathbb{R}_{>0})$, being $\tau$-noncausal with respect to $\Omega(\mathbb{R}_{\geq 0}, U)$, in such a way that*

$$\begin{aligned} &e'P(t)A(t,q,y(t))e + \tfrac{1}{2}e'\dot{P}(t)e \leq \phi(t)|H(t,y(t))e|^2 - \bar{d}(t)e'P(t)e, \\ &\forall t \geq \bar{t}_0,\, e \in \mathbb{R}^n,\, q \in Q(t), \text{provided that (2.10) holds} \end{aligned} \tag{2.12}$$

**Proof.** Let $t_0$, $\bar{t}_0$, $\tau_0$ and $\tau$ as given in our statement and let $y \in \Omega(\mathbb{R}_{\geq 0}, U) \cap C^0([t_0, \infty); \mathbb{R}^k)$. To simplify the proof, we may consider two cases.

**Case I**: (2.10) holds; i.e. $y \in \Omega(\mathbb{R}_{\geq 0}, U \cap W) \cap C^0([t_0, \infty); \mathbb{R}^k)$.

For any $y(\cdot)$ satisfying (2.10) and for each $t \geq \bar{t}_0$, $q \in \mathbb{R}^\ell$ and $e \in \mathbb{R}^n$ we define:

$$D_y(t, q, e) := e'P(t)A(t,q,y(t))e + \tfrac{1}{2}e'\dot{P}(t)e + \bar{d}(t)e'P(t)e\,; \tag{2.13a}$$

$$K(t) := \{w \in \mathbb{R}^n : |w| = 1, D_y(t, q, w) < 0, \forall q \in Q(t)\} \tag{2.13b}$$

For those $t \geq \bar{t}_0$ for which $\operatorname{rank} H(t, y(t)) < n$, the set $K(t)$ is nonempty, since it includes all vectors $w \in \mathbb{R}^n : |w| = 1$ with $w \in \ker H(t, y(t))$. Indeed, let $w \in \ker H(t, y(t))$ for certain nonzero $w \in \mathbb{R}^n$. Then, then by using (2.9),(2.10),(2.11),(2.13a) and by taking into account that $P(\cdot)$ is positive definite, we deduce $D_y(t, q, w) \leq (\bar{d}(t) - d(t))w'P(t)w < 0$, for all $q \in Q(t)$, therefore $w \in K(t)$, which establishes that $K(t) \neq \emptyset$. The above discussion asserts that:

$$(w \in \ker H(t, y(t)) \text{ and } |w|=1) \Rightarrow w \in K(t) \tag{2.14}$$

In the sequel, for each $t \geq \bar{t}_0$ we adopt the notation $K^c(t)$ to denote the complement of $K(t)$ with respect to the unit sphere centered at $0 \in \mathbb{R}^n$, namely, $K^c(t) := \{w \in \mathbb{R}^n : |w|=1, w \notin K(t)\}$, which by virtue of (2.13b) is written:

$$K^c(t) = \{w \in \mathbb{R}^n : |w|=1 \text{ with } D_y(t, q, w) \geq 0, \text{ for some } q \in Q(t)\} \tag{2.15}$$

Notice that $K^c(t)$ is empty, if and only if

$$D_y(t, q, w) < 0, \text{ for every } q \in Q(t) \text{ and } w \in \mathbb{R}^n \text{ with } |w|=1 \tag{2.16}$$

We now prove that for every $t \geq \bar{t}_0$ the set $K^c(t)$ is closed. Let $t \geq \bar{t}_0$ and without any loss of generality let us assume that $K^c(t) \neq \emptyset$. We prove that for every sequence $(w_\nu)_{\nu \in \mathbb{N}} \subset K^c(t)$ with $w_\nu \to w$ we have $w \in K^c(t)$. Indeed, since $w_\nu \in K^c(t)$, (2.15) asserts that there exists $q_\nu \in Q(t)$ with $D_y(t, q_\nu, w_\nu) \geq 0$. By recalling the CP for the map $Q(\cdot)$ we may assume that $q_\nu \to q$ for certain $q \in Q(t)$. Since $(q_\nu, w_\nu) \to (q, w)$, continuity of $D_y(t, \cdot, \cdot)$ implies $D_y(t, q_\nu, w_\nu) \to D_y(t, q, w) \geq 0$, thus by (2.15) $w \in K^c(t)$ and consequently $K^c(t)$ is closed. Next, consider the map $\omega : [\bar{t}_0, \infty) \to [0, \infty]$ defined as

$$\omega(t) := \begin{cases} \min\{|H(t, y(t))w| : w \in K^c(t)\}, & \text{if } K^c(t) \neq \emptyset \\ \infty & \text{if } K^c(t) = \emptyset \end{cases} \tag{2.17}$$

Notice that the set $\{|H(t, y(t))w| : w \in K^c(t)\}$ is compact, whenever $K^c(t) \neq \emptyset$, thus $\omega(\cdot)$ is well defined, $\tau_0$-noncausal with respect to $\Omega(\mathbb{R}_{\geq 0}, U)$ and, due to (2.14), satisfies $\omega(t) > 0$ for every $t \geq \bar{t}_0$. Moreover, we show that for every $T > \bar{t}_0$ it holds:

$$\inf\{\omega(t) : t \in [\bar{t}_0, T]\} > 0 \tag{2.18}$$

Indeed, suppose on the contrary that $\omega(t_\nu) \to 0$ for some $t_\nu \to t \in [\bar{t}_0, T]$. Then by taking into account (2.17) we may assume that without any loss of generality it holds $K^c(t_\nu) \neq \emptyset$ for every $\nu \in \mathbb{N}$, hence, there exist $w_\nu \in K^c(t_\nu)$ such that $|H(t_\nu, y(t_\nu))w_\nu| \to 0$. Since $|w_\nu|=1$, we may also assume that, without any loss of generality, there exists $w \in \mathbb{R}^n$ with $|w|=1$ and $w_\nu \to w$, therefore, continuity of $H(\cdot, y(\cdot))$ implies $H(t, y(t))w = 0$. It follows by virtue of (2.14) that $w \in K(t)$. On the other hand, $w_\nu \in K^c(t_\nu)$, hence $D_y(t_\nu, q_\nu, w_\nu) \geq 0$ for some $q_\nu \in Q(t_\nu)$. The latter in conjunction with the compactness property for $Q(\cdot)$ and continuity of $D_y(\cdot, \cdot, \cdot)$ implies that $D_y(t_\nu, q_\nu, w_\nu) \to D_y(t, q, w) \geq 0$ for some $q \in Q(t)$, therefore $w \in K^c(t)$, which is a contradiction. We conclude that (2.18) is fulfilled. Now, let $\bar{\omega} : [\bar{t}_0, \infty) \to [0, \infty)$ defined as

$$\bar{\omega}(t) := \begin{cases} \dfrac{1}{\omega^2(t)}, & K^c(t) \neq \emptyset \\ 0, & K^c(t) = \emptyset \end{cases} \tag{2.19}$$

The function $\bar{\omega}(\cdot)$ is $\tau_0$-noncausal with respect to $\Omega(\mathbb{R}_{\geq 0}, U)$ and by (2.18) and (2.19) it follows that for any $T > \bar{t}_0$ there exists a constant $M := M(T) > 0$ such that

$$\sup\{\bar{\omega}(t) : t \in [\bar{t}_0, T]\} \leq M \tag{2.20}$$

Define:

$$C(t) := \sup\{\overline{\omega}(t)(|P(t)\| A(t,q,y(t))| + \tfrac{1}{2}|\dot{P}(t)| + |\overline{d}(t)\| P(t)|) : q \in Q(t)\} \quad (2.21)$$

The map $C(\cdot)$ is $\tau_0$-noncausal with respect to $\Omega(\mathbb{R}_{\geq 0}, U)$. Furthermore, by taking into account (2.20) and (2.21) and recalling again the CP for $Q(\cdot)$, it follows that for any $T > \overline{t_0}$ the set $\bigcup_{t \in [\overline{t_0}, T]} Q(t)$ is bounded and there exists a constant $\overline{M} := \overline{M}(T) > 0$ such that

$$\sup\{C(t) : t \in [\overline{t_0}, T]\} \leq \overline{M} \quad (2.22)$$

Since $\tau > \tau_0$, we can apply the same arguments with those in the proof of Proposition 2.1(i) in [49], in order to build a function $\phi \in C^1([\overline{t_0}, \infty); \mathbb{R}_{>0})$, which is $\tau$-noncausal with respect to $\Omega(\mathbb{R}_{\geq 0}, U)$ and satisfies:

$$\phi(t) > C(t), \ \forall t \geq \overline{t_0} \quad (2.23)$$

Notice that, due to (2.10), the desired (2.12) is equivalent to

$$w'P(t)A(t,q,y(t))w + \tfrac{1}{2}w'\dot{P}(t)w \leq \phi(t)|H(t,y(t))w|^2 - \overline{d}(t)w'P(t)w \quad (2.24)$$
$$\forall t \geq \overline{t_0}, w \in \mathbb{R}^n \text{ with } |w| = 1, q \in Q(t)$$

Thus, in order to establish (2.12), it suffices to show that (2.24) holds. We distinguish two cases. First, consider those $t \geq \overline{t_0}$ for which $K^c(t) \neq \emptyset$. Then, if $K(t) \neq \emptyset$ and $w \in K(t)$, (2.24) is a consequence of (2.13). When $w \in K^c(t)$, in order to show (2.24), it suffices by virtue of (2.17) to prove:

$$\sup\{|P(t)\| A(t,q,y(t))| + \tfrac{1}{2}|\dot{P}(t)| + |\overline{d}(t)\| P(t)| : q \in Q(t)\} \leq \phi(t)\omega^2(t) \quad (2.25)$$

The latter is a consequence of (2.19),(2.21) and (2.23). Now, for those $t \geq \overline{t_0}$ for which $K^c(t) = \emptyset$ it follows by (2.16) that $D_y(t,q,w) < 0$, for every $w \in \mathbb{R}^n$ with $|w| = 1$ and $q \in Q(t)$, thus (2.24) is a consequence of (2.13a).

**Case II:** $y \in \Omega(\mathbb{R}_{\geq 0}, U) \cap C^0([t_0, \infty); \mathbb{R}^k)$.

For the general case above, we make some appropriate modifications to the procedure used in the Case I. Particularly, define $D_y(\cdot, \cdot, \cdot,)$, $K(\cdot)$ and $\omega(\cdot)$ as in (2.13a),(2.13b) and (2.17), respectively and notice that $\omega(t) \geq 0$ for all $t \geq \overline{t_0}$. Instead of (2.19), we define $\overline{\omega} : [\overline{t_0}, \infty) \to [0, \infty)$ as follows:

$$\overline{\omega}(t) := \begin{cases} \dfrac{1}{\omega^2(t)}, & \text{provided that } \inf\{\omega(s) : s \in [t - \tfrac{\tau - \tau_0}{2}, t + \tfrac{\tau - \tau_0}{2}] \cap [\overline{t_0}, \infty)\} > 0 \\ 0, & \text{otherwise} \end{cases} \quad (2.26)$$

Obviously, $\overline{\omega}(\cdot)$ is $(\tau_0 + \tfrac{\tau - \tau_0}{2})$-noncausal with respect to $\Omega(\mathbb{R}_{\geq 0}, U)$. Also, consider the function $C(\cdot)$ as precisely defined by (2.21). Since $\tau > \tau_0 + \tfrac{\tau - \tau_0}{2}$ we can find a function $\phi \in C^1([\overline{t_0}, \infty); \mathbb{R}_{>0})$, being $\tau$-noncausal with respect to $\Omega(\mathbb{R}_{\geq 0}, U)$, in such a way that (2.23) is satisfied. Then, when (2.10) holds, it follows by taking into account (2.18) that the function $\overline{\omega}(\cdot)$ as defined by (2.19) coincides with $\overline{\omega}(\cdot)$ as the latter is given by (2.26), hence, by repeating the same arguments used in the proof of Case I, it follows that (2.12) is fulfilled. ∎

**Corollary 2.1.** *Consider the pair* $(H, A)$ *as given in (2.8a,b) with* $p = k$, $H(t, y) := H(t)$ *and* $A(t, q, y)$ *as involved in hypotheses A1 and A2. Suppose that A2 is fulfilled and consider the constants* $R$, $\varepsilon_R$, $\xi$, $\overline{t_0} \geq t_0 \geq 0$ *and* $\tau_0 > 0$ *and the mappings* $Q_R := Q_{R, t_0, \xi}$, $y \in Y(\mathbb{R}_{\geq 0}, M) \cap C^0([t_0, \infty); \mathbb{R}^k)$, $P_R(\cdot) = P_{R,t_0,\bar{t}_0,\tau_0,\xi,y}(\cdot)$ and $d_R(\cdot) = d_{R,t_0,\bar{t}_0,\tau_0,\xi,y}(\cdot)$, as precisely determined in A2. Then for every $\bar{\varepsilon}_R > \varepsilon_R$, $\tau > \tau_0$ there exist functions $\bar{d}_R \in C^0([\bar{t}_0, \infty); \mathbb{R})$ with

$$\bar{d}_R(t) < d_R(t), \ \forall t \geq \bar{t}_0 ; \tag{2.27a}$$

$$\int_{\bar{t}_0}^{t} \bar{d}_R(s)ds > -\bar{\varepsilon}_R, \ \forall t \geq \bar{t}_0 ; \ \int_{\bar{t}_0}^{\infty} \bar{d}_R(s)ds = \infty \tag{2.27b}$$

and $\phi_R \in C^1([\bar{t}_0, \infty); \mathbb{R}_{>0})$ such that

$$e'P_R(t)A(t,q,y(t))e + \tfrac{1}{2}e'\dot{P}_R(t)e - \phi_R(t)|H(t)e|^2 \leq -\bar{d}_R(t)e'P_R(t)e,$$
$$\forall t \geq \bar{t}_0, \ e \in \mathbb{R}^n, \ q \in Q_R(t), \text{provided that (2.7) holds} \tag{2.27c}$$

both being $\tau$-noncausal with respect to $Y(\mathbb{R}_{\geq 0}, M)$. Particularly, the function $\phi_R$ is any $\tau$-noncausal $C^1([\bar{t}_0, \infty); \mathbb{R}_{>0})$ map satisfying

$$\phi_R(t) > C_R(t), \ \forall t \geq \bar{t}_0 ; \tag{2.28a}$$

$$C_R(t) := \sup\{\bar{\omega}_R(t)(|P_R(t)\| A(t,q,y(t))| + \tfrac{1}{2}|\dot{P}_R(t)| + |\bar{d}_R(t)\| P_R(t)|) : q \in Q_R(t)\} ; \tag{2.28b}$$

$$\bar{\omega}_R(t) := \begin{cases} \dfrac{1}{\omega_R^2(t)}, & \text{provided that } \inf\{\omega_R(s) : s \in [t-\tfrac{\tau-\tau_0}{2}, t+\tfrac{\tau-\tau_0}{2}] \cap [\bar{t}_0, \infty)\} > 0 \\ 0, & \text{otherwise} \end{cases} ; \tag{2.28c}$$

$$\omega_R(t) := \begin{cases} \min\{|H(t,y(t))w| : w \in K_R^c(t)\}, & \text{if } K_R^c(t) \neq \emptyset \\ \infty & \text{if } K_R^c(t) = \emptyset \end{cases} ; \tag{2.28d}$$

$$K_R(t) := \{w \in \mathbb{R}^n : |w| = 1, w'P_R(t)A(t,q,y(t))w + \tfrac{1}{2}w'\dot{P}_R(t)w + \bar{d}_R(t)w'P_R(t)w \leq 0, \forall q \in Q_R(t)\},$$
$$t \geq \bar{t}_0. \tag{2.28e}$$

**Proof.** Let $\bar{\varepsilon}_R > \varepsilon_R$ and define

$$\bar{d}_R(t) := d_R(t) - \frac{2(\bar{\varepsilon}_R - \varepsilon_R)}{\pi(1+t^2)}, \ t \geq \bar{t}_0 \tag{2.29}$$

Then by taking into account (2.6b) and (2.29) it follows that (2.27a,b) hold. Also, notice that $\bar{d}_R \in C^0([\bar{t}_0, \infty); \mathbb{R})$ is $\tau_0$-noncausal with respect to $Y(\mathbb{R}_{\geq 0}, M)$. Furthermore, due to (2.6c) and (2.7), property H1 holds for the pair $(H, A)$ with $U := M$, $W := B_R$, $\Omega(\mathbb{R}_{\geq 0}, U) := Y(\mathbb{R}_{\geq 0}, M)$, $P := P_R$, $d := d_R$ and $Q := Q_R$. It follows, by taking into account (2.27a) and Proposition 2.1, that for the given $\bar{t}_0 \geq t_0 \geq 0$, $\tau > \tau_0 > 0$ and $y \in Y(\mathbb{R}_{\geq 0}, M) \cap C^0([t_0, \infty); \mathbb{R}^k)$, any $\tau$-noncausal with respect to $Y(\mathbb{R}_{\geq 0}, M)$ function $\phi_R \in C^1([\bar{t}_0, \infty); \mathbb{R}_{>0})$ as defined in (2.28a)-(2.28e), satisfies the desired (2.27c). ∎

The following proposition partially generalizes Proposition 2.1 in [49] providing sufficient conditions for the solvability of the observer design problem for the case (2.1) under the additional hypothesis that all initial states belong to a given known compact set. Specifically, it is a priori known that the initial states $x_0$ of (2.1a) belong to the compact ball $B_R$ of radius $R > 0$ centered at zero.

**Proposition 2.2.** *Consider the system (2.1) and let $M$ be a nonempty subset of $\mathbb{R}^n$ with $0 \in \mathrm{cl}\, M$ such that system (2.1a) is M-forward complete. For the initial state $x_0 \in M$ of (2.1a), assume that $|x_0| \leq R$ for some known constant $R > 0$ and assume that properties A1 and A2 hold with $R$ as above and for certain constant $L > 1$. Then the AC-ODP is solvable for (2.1) with respect to $Y(\mathbb{R}_{\geq 0}, B_R \cap M)$. Particularly, for each $\bar{t}_0 \geq t_0 \geq 0$, $\tau > 0$, $\bar{\varepsilon}_R > \varepsilon_R$ ($\varepsilon_R$ being the constant involved in A2) and constant $\xi$ satisfying*

$$\xi \geq \beta(\bar{t}_0, t_0, R)\sqrt{L} \exp[\bar{\varepsilon}_R] \tag{2.30}$$

*and for every* $y \in Y(\mathbb{R}_{\geq 0}, B_R \cap M) \cap C^0([t_0, \infty); \mathbb{R}^k)$ *there exist a time-varying symmetric matrix* $P_R \in C^1([\bar{t}_0, \infty); \mathbb{R}^{n \times n})$ *satisfying (2.6a) and functions* $\bar{d}_R \in C^0([\bar{t}_0, \infty); \mathbb{R})$, $\phi_R \in C^1([\bar{t}_0, \infty); \mathbb{R}_{>0})$, *being* $\tau$-*noncausal with respect to* $Y(\mathbb{R}_{\geq 0}, B_R \cap M)$, *in such a way that (2.27) holds with* $Q_R := Q_{R,t_0,\xi}$ *as precisely given in A1. It turns out that the AC-ODP is solvable for (2.1) with respect to* $Y(\mathbb{R}_{\geq 0}, B_R \cap M)$. *Particularly:*

(i) *The system*

$$\dot{z} = F(t, z, y) + \phi_R(t) P_R^{-1}(t) H'(t)(y - H(t)z),\quad (2.31a)$$
$$\text{with initial } z(\bar{t}_0) = 0 \quad (2.31b)$$

*is forward complete and is an observer for (2.1);*

(ii) *The error between the trajectory* $x(\cdot) := x(\cdot, t_0, x_0)$ *of (2.1a) and the trajectory* $z(\cdot) := z(\cdot, \bar{t}_0, 0)$ *of the observer (2.31) is given by*

$$|e(t)| < \beta(\bar{t}_0, t_0, R)\sqrt{L}\exp\left[\int_{\bar{t}_0}^t -\bar{d}_R(s)ds\right],\ \forall t \geq \bar{t}_0 \quad (2.32)$$

**Proof.** Let $\bar{t}_0 \geq t_0 \geq 0$, $\tau > 0$, $\bar{\varepsilon}_R > \varepsilon_R$ and let $\xi$ be a constant satisfying (2.30). Also, consider $y(\cdot)$ satisfying (2.7), namely $y \in Y(\mathbb{R}_{\geq 0}, B_R \cap M) \cap C^0([t_0, \infty); \mathbb{R}^k)$. Then according to A2 there exist a time-varying symmetric matrix $P_R \in C^1([\bar{t}_0, \infty); \mathbb{R}^{n \times n})$ and a function $d_R \in C^0([\bar{t}_0, \infty); \mathbb{R})$, both $\tau_0$-noncausal with respect to $Y(\mathbb{R}_{\geq 0}, M)$ for certain $\tau_0 \in (0, \tau)$, in such a way that (2.6a,b,c) hold. Therefore, since $\bar{\varepsilon}_R > \varepsilon_R$ and $\tau > \tau_0$, Corollary 2.1 asserts that there exist $\bar{d}_R \in C^0([\bar{t}_0, \infty); \mathbb{R})$ and $\phi_R \in C^1([\bar{t}_0, \infty); \mathbb{R}_{>0})$ both $\tau$-noncausal with respect to $Y(\mathbb{R}_{\geq 0}, M)$ - hence, by (2.7) both being $\tau$-noncausal with respect to $Y(\mathbb{R}_{\geq 0}, B_R \cap M)$ - in such a way that (2.27a,b,c) hold. Let $e := x - z$. Then by invoking (2.1b) the error equation between (2.31a) and (2.1a) is written:

$$\begin{aligned}\dot{e} &= F(t, x, y) - F(t, z, y) - \phi_R(t)P_R^{-1}(t)H'(t)(y - H(t)z)\\ &= F(t, x, y) - F(t, z, y) - \phi_R(t)P_R^{-1}(t)H'(t)(H(t)x - H(t)z)\\ &= F(t, x, y) - F(t, z, y) - \phi_R(t)P_R^{-1}(t)H'(t)H(t)e,\ t \geq \bar{t}_0 \end{aligned} \quad (2.33)$$

Denote by $T_{\max}$ the maximum time for which the solution $e(\cdot)$ of (2.33) with initial condition

$$e(\bar{t}_0) = x(\bar{t}_0, t_0, x_0) - z(\bar{t}_0) \overset{(2.31b)}{=} x(\bar{t}_0, t_0, x_0) \quad (2.34)$$

exists on the interval $[\bar{t}_0, T_{\max})$. Notice that, due to (1.4),(2.30) and (2.34), it holds that $|e(\bar{t}_0)| < \xi$. We claim that $|e(t)| < \xi$, for every $t \in [\bar{t}_0, T_{\max})$ and therefore $T_{\max} = \infty$. Indeed, suppose on the contrary that there exists a time $\bar{t} \in (\bar{t}_0, T_{\max})$ such that

$$|e(\bar{t})| = \xi; \quad (2.35a)$$
$$|e(t)| < \xi,\ \forall t \in [\bar{t}_0, \bar{t}) \quad (2.35b)$$

Now, by invoking (1.4) and (2.4) we have:

$$y(t) \in Y_R(t),\ \forall t \in [\bar{t}_0, \bar{t}] \quad (2.36)$$

and, by recalling (2.5) in A1, it follows from (2.35a,b) and (2.36) that

$$\begin{aligned}\Delta F(t, x(t), z(t); y(t)) &= A(t, q, y(t))(x(t) - z(t)), \forall t \in [\bar{t}_0, \bar{t}]\\ &\text{and for some } q := q(t) \in Q_R(t)\end{aligned} \quad (2.37)$$

thus, by evaluating the time-derivative $\dot{V}$ of $V(t,e) := \frac{1}{2} e' P_R(t) e$, $e \in \mathbb{R}^n$ along the trajectories $e(\cdot) = e(\cdot, \bar{t}_0, e(\bar{t}_0), x(\cdot))$ of (2.33) and by exploiting (2.37), we find:

$$\dot{V} = \tfrac{1}{2} e'(t) \dot{P}_R(t) e(t) + e'(t) P_R(t) A(t,q,y(t)) e(t) - \phi_R(t) |H(t) e(t)|^2, \quad \forall t \in [\bar{t}_0, \bar{t}] \tag{2.38}$$

By (2.27c) and (2.38) it follows that

$$\dot{V} \le -2 \bar{d}_R(t) V, \quad \forall t \in [\bar{t}_0, \bar{t}] \tag{2.39}$$

Consequently, from (2.39) we deduce $\tfrac{1}{2} e'(t) P_R(t) e(t) \le \tfrac{1}{2} e'(\bar{t}_0) P_R(\bar{t}_0) e(\bar{t}_0) \exp\left[\int_{\bar{t}_0}^{t} -2 \bar{d}_R(s) ds\right]$, therefore, by taking into account the first inequality of (2.6a), we have:

$$|e(t)| \le |e(\bar{t}_0)| \, \|P_R(\bar{t}_0)\|^{\frac{1}{2}} \exp\left[\int_{\bar{t}_0}^{t} -\bar{d}_R(s) ds\right], \quad \forall t \in [\bar{t}_0, \bar{t}] \tag{2.40}$$

Notice that, since $z(\bar{t}_0) = 0$, (1.4) implies:

$$|e(\bar{t}_0)| = |x(\bar{t}_0, t_0, x_0)| \le \beta(\bar{t}_0, t_0, R) \tag{2.41}$$

By taking into account the second inequality of (2.6a), as well as (2.27b),(2.30),(2.40) and (2.41), we get:

$$|e(t)| < \beta(\bar{t}_0, t_0, R) \sqrt{L} \exp[\bar{\varepsilon}_R] = \xi, \quad \forall t \in [\bar{t}_0, \bar{t}] \tag{2.42}$$

thus $|e(\bar{t})| < \xi$, which contradicts (2.35a). It follows that the solution $e(\cdot) = e(\cdot, t_0, e(t_0), x(\cdot))$ of (2.33) satisfies $|e(t)| < \xi$ for every $t \in [\bar{t}_0, T_{\max})$ and therefore $T_{\max} = \infty$. Finally, by recalling the second inequality of (2.6a),(2.40) and (2.41) it also follows that (2.32) is fulfilled, which, by virtue of the equality in (2.27b), asserts that (1.8) holds with $M := B_R \cap M$. We conclude that the AC-ODP is solvable with respect to $Y(\mathbb{R}_{\ge 0}, B_R \cap M)$ and (2.31) is an observer for (2.1). ∎

We now generalize the result of Proposition 2.2 by establishing sufficient conditions for the existence of a switching observer exhibiting the state determination of (2.1), without any priori information concerning the initial condition. We make the following hypothesis:

**A3** *Assume that there exist an integer $\ell \in \mathbb{N}$, a map $A(\cdot, \cdot, \cdot)$ as in (2.2) and a constant $L > 1$, in such a way that for every $R > 0$ both hypotheses A1 and A2 hold.*

We now state our main result for the solvability of the AC-SODP.

**Proposition 2.3.** *In addition to hypothesis of M-forward completeness for (2.1a), assume that system (2.1) satisfies A3. Then the AC-SODP is solvable for (2.1) with respect to $Y(\mathbb{R}_{\ge 0}, M)$.*

**Proof.** Consider the system (2.1) initiated at $t_0 \ge 0$. For arbitrary $\tau > 0$ we proceed to the construction of a sequence of times $(t_m)_{m \in \mathbb{N}}$ and appropriate sequences of mappings $(g_m)_{m \in \mathbb{N}}$, being $\tau$-noncausal with respect to $Y(\mathbb{R}_{\ge 0}, M)$, and satisfying the requirements of Definition 1.3. Let $\tau > 0$ and let $y \in Y(\mathbb{R}_{\ge 0}, M) \cap C^0([t_0, \infty); \mathbb{R}^k)$ be the output of (2.1). Also, let $L > 1$ be a constant such that for every $R > 0$ properties A1 and A2 hold, according to our Assumption A3.

**Claim:** For $L$, $\tau$ and $y(\cdot)$ as above and for any $m \in \mathbb{N}$, there exist positive constants $\varepsilon_i$, $\xi_i$ and $t_i$, set-valued mappings $Q_i(\cdot)$ satisfying the CP, positive definite mappings $P_i \in C^1([t_{i-1}, \infty); \mathbb{R}^{n \times n})$ and functions $\bar{d}_i \in C^0([t_{i-1}, \infty); \mathbb{R})$ and $\phi_i \in C^1([t_{i-1}, \infty); \mathbb{R}_{>0})$, $i = 1, 2, \ldots, m$, all $\tau$-noncausal with respect to $Y(\mathbb{R}_{\ge 0}, M)$, in such a way that

$$\xi_i := \beta(t_{i-1}, t_0, i)\sqrt{L}\exp[\overline{\varepsilon}_i]; \quad \overline{\varepsilon}_i := 2\varepsilon_i; \tag{2.43}$$

$$\int_{t_{i-1}}^{t} \overline{d}_i(s)ds > -\overline{\varepsilon}_i, \quad \forall t \geq t_{i-1}; \quad \int_{t_{i-1}}^{\infty} \overline{d}_i(s)ds = \infty; \tag{2.44}$$

$$t_1 := t_0 \text{ for the case } m = 1;$$

$$t_i := \min\left\{T \geq t_{i-1} + 1 : \exp\left\{\int_{t_{i-1}}^{t} -\overline{d}_i(s)ds\right\} \leq \frac{1}{\beta(t_{i-1}, t_0, i)i\sqrt{L}}, \text{ for all } t \geq T\right\},$$

$$i = 2, \ldots, m \text{ for the case } m \geq 2$$
$$\tag{2.45}$$

and further

$$e'P_i(t)A(t, q, y(t))e + \tfrac{1}{2}e'\dot{P}_i(t)e - \phi_i(t)\,|\,H(t)e\,|^2 \leq -\overline{d}_i(t)e'P_i(t)e,$$
$$\forall t \geq t_{i-1}, e \in \mathbb{R}^n, q \in Q_i(t), \text{provided that (2.7) holds with } R = i \tag{2.46}$$

It turns out by (2.45) that the sequence $\{t_i\}$ above satisfies:

$$\lim_{i \to \infty} t_i = \infty \tag{2.47}$$

In order to establish our claim, we proceed by induction as follows. We first consider the case $m := 1$. Due to Assumption A3, we may apply property A2 with $R := 1$, which asserts that there exists a constant $\varepsilon_1 > 0$ such that, if for the given $t_0$ above we define:

$$\xi_1 := \beta(t_0, t_0, 1)\sqrt{L}\exp[\overline{\varepsilon}_1]; \quad \overline{\varepsilon}_1 := 2\varepsilon_1, \tag{2.48}$$

then there exists a set-valued map $Q_1 := Q_{1, t_0, \xi_1}$ satisfying the CP in such a way that for $\overline{t}_0 := t_0 \geq t_0$, $\tau_0 := \tfrac{\tau}{2}$ and output $y(\cdot)$ as above, a pair of mappings

$$d_1 \in C^0([t_0, \infty); \mathbb{R}) \text{ and } P_1 \in C^1([t_0, \infty); \mathbb{R}^{n \times n})$$

can be found, both $\tau_0$-noncausal with respect to $Y(\mathbb{R}_{\geq 0}, M)$ and $P_1(\cdot)$ being positive definite, such that for $R = 1$ conditions (2.6a,b) hold and further (2.6c) is fulfilled provided that (2.7) holds. Then, by invoking Corollary 2.1 with $R = 1$, for the given $\overline{t}_0 \geq t_0$, $\tau > \tau_0 (= \tfrac{\tau}{2})$, $\xi_1$, $\overline{\varepsilon}_1$ and $y(\cdot)$ as above, there exist functions

$$\overline{d}_1 \in C^0([t_0, \infty); \mathbb{R}) \text{ and } \phi_1 \in C^1([t_0, \infty); \mathbb{R}_{>0}) \tag{2.49}$$

both being $\tau$-noncausal with respect to $Y(\mathbb{R}_{\geq 0}, M)$, which satisfy:

$$\int_{t_0}^{t} \overline{d}_1(s)ds > -\overline{\varepsilon}_1, \quad \forall t \geq t_0; \quad \int_{t_0}^{\infty} \overline{d}_1(s)ds = \infty \tag{2.50a}$$

and further

$$e'P_1(t)A(t, q, y(t))e + \tfrac{1}{2}e'\dot{P}_1(t)e - \phi_1(t)\,|\,H(t)e\,|^2 \leq -\overline{d}_1(t)e'P_1(t)e,$$
$$\forall t \geq t_0, e \in \mathbb{R}^n, q \in Q_1(t), \text{provided that (2.7) holds with } R = 1 \tag{2.50b}$$

Suppose next that (2.43)-(2.46) are fulfilled for certain $m \in \mathbb{N}$. For $m := m + 1$ we define the corresponding constants $\varepsilon_{m+1}$, $\xi_{m+1}$, $t_{m+1}$ and mappings $Q_{m+1}(\cdot)$, $P_{m+1}(\cdot)$, $\overline{d}_{m+1}(\cdot)$ and $\phi_{m+1}(\cdot)$ as follows. Due to Assumption A3, we may apply again A2 with $R := m + 1$, which asserts existence of a constant $\varepsilon_{m+1} > 0$ such that, if we define:

$$\xi_{m+1} := \beta(t_m, t_0, m+1)\sqrt{L}\exp[\overline{\varepsilon}_{m+1}]; \quad \overline{\varepsilon}_{m+1} := 2\varepsilon_{m+1}, \tag{2.51}$$

then a set-valued map $Q_{m+1} := Q_{m+1,t_0,\xi_{m+1}}$ satisfying the CP can be found in such a way that for $\bar{t}_0 := t_m \geq t_0$, $\tau_0 := \frac{\tau}{2}$ and output $y(\cdot)$ as above there exists a pair of $\tau_0$-noncausal with respect to $Y(\mathbb{R}_{\geq 0}, M)$ mappings

$$d_{m+1} \in C^0([t_m, \infty); \mathbb{R}) \text{ and } P_{m+1} \in C^1([t_m, \infty); \mathbb{R}^{n \times n})$$

$P_{m+1}(\cdot)$ being positive definite, such that for $R = m+1$ conditions (2.6a,b) hold and further (2.6c) is fulfilled, provided that (2.7) holds. Then, by invoking again Corollary 2.1 with $R = m+1$, for the given $\bar{t}_0 \geq t_0$, $\tau > \tau_0 (= \frac{\tau}{2})$, $\xi_{m+1}$, $\bar{\varepsilon}_{m+1}$ and $y(\cdot)$, there exist functions

$$\bar{d}_{m+1} \in C^0([t_m, \infty); \mathbb{R}) \text{ and } \phi_{m+1} \in C^1([t_m, \infty); \mathbb{R}_{>0}) \tag{2.52}$$

both being $\tau$-noncausal with respect to $Y(\mathbb{R}_{\geq 0}, M)$, which satisfy:

$$\int_{t_m}^t \bar{d}_{m+1}(s) ds > -\bar{\varepsilon}_{m+1}, \ \forall t \geq t_m; \ \int_{t_m}^\infty \bar{d}_{m+1}(s) ds = \infty \tag{2.53a}$$

and further

$$e' P_{m+1}(t) A(t, q, y(t)) e + \tfrac{1}{2} e' \dot{P}_{m+1}(t) e - \phi_{m+1}(t) |H(t)e|^2 \leq -\bar{d}_{m+1}(t) e' P_{m+1}(t) e, \tag{2.53b}$$
$$\forall t \geq t_m, e \in \mathbb{R}^n, q \in Q_{m+1}(t), \text{ provided that (2.7) holds with } R = m+1$$

Finally, due to (2.53a), we may define:

$$t_{m+1} := \min\left\{ T \geq t_m + 1 : \exp\left\{ \int_{t_m}^t -\bar{d}_{m+1}(s) ds \right\} \leq \frac{1}{\beta(t_m, t_0, m+1)(m+1)\sqrt{L}}, \text{ for all } t \geq T \right\} \tag{2.54}$$

and this completes the establishment of our claim.

We are now in a position, by taking into account (2.43)-(2.47), to construct the desired switching observer exhibiting the state determination of (1.8), according to the Definition 1.3. Consider for each $m \in \mathbb{N}$ the system

$$\dot{z}_m = g_m(t, z_m, y(t)), \ t \in [t_{m-1}, t_{m+1}], \tag{2.55a}$$
$$\text{with initial } z_m(t_{m-1}) = 0 \tag{2.55b}$$

with dynamics:

$$g_m(t, z, y) \begin{cases} := F(t, z, y) + \phi_m(t) P_m^{-1}(t) H'(t)(y - H(t)z) \\ \qquad \text{for } |z| \leq \zeta_m, \ t \in [t_{m-1}, t_{m+1}], \\ := (F(t, z, y) + \phi_m(t) P_m^{-1}(t) H'(t)(y - H(t)z)) \dfrac{2\zeta_m - |z|}{\zeta_m} \\ \qquad \text{for } \zeta_m \leq |z| \leq 2\zeta_m, \ t \in [t_{m-1}, t_{m+1}], \\ := 0, \qquad \text{for } |z| \geq 2\zeta_m, \ t \in [t_{m-1}, t_{m+1}] \end{cases} \tag{2.56a}$$

$$\zeta_m := \beta(t_{m+1}, t_0, m) + \xi_m \tag{2.56b}$$

where $(t_m)_{m \in \mathbb{N}}$, $(\xi_m)_{m \in \mathbb{N}}$, $(\phi_m)_{m \in \mathbb{N}}$ and $(P_m)_{m \in \mathbb{N}}$ are determined in (2.43)-(2.46). Notice that $g_m(\cdot)$ is $C^0([t_{m-1}, t_{m+1}] \times \mathbb{R}^n \times \mathbb{R}^k; \mathbb{R}^n)$ and further it is locally Lipschitz on $(z, y) \in \mathbb{R}^n \times \mathbb{R}^k$, since it is the product of the following locally Lipschitz on $(z, y) \in \mathbb{R}^n \times \mathbb{R}^k$ functions:

$$f_{1,m}(t, z, y) := F(t, z, y) + \phi_m(t) P_m^{-1}(t) H'(t)(y - H(t)z)$$

$$f_{2,m}(z,y) := \begin{cases} 1, & |z| \leq \zeta_m \\ \dfrac{2\zeta_m - |z|}{\zeta_m}, & \zeta_m \leq |z| \leq 2\zeta_m \\ 0, & |z| \geq 2\zeta_m \end{cases}$$

Thus for each $m \in \mathbb{N}$ the system (2.55a) with dynamics (2.56) has a unique solution $z_m(\cdot)$ with $z_m(t_{m-1}) = 0$, which exists for $t \geq t_{m-1}$ near $t_{m-1}$. We show that $z_m(\cdot)$ is defined on the interval $[t_{m-1}, t_{m+1}]$ and particularly, we claim that $|z_m(t)| < 2\zeta_m$ for all $t \in [t_{m-1}, t_{m+1}]$. Indeed, otherwise, since $z(t_{m-1}) = 0$, there exists $\tilde{t} \in (t_{m-1}, t_{m+1}]$ such that $|z_m(\tilde{t})| = 2\zeta_m$, whereas $|z_m(t)| < 2\zeta_m$ for $t \in (t_{m-1}, \tilde{t})$. Then by (2.55a) and (2.56a) it follows that $\dot{z}_m(t) = 0$, hence $|z_m(t)| = 2\zeta_m$ for each $t$ around $\tilde{t}$, which is a contradiction.

We are now in a position to establish the desired (1.8). Let $Z : [t_0, \infty) \to \mathbb{R}^n$ as defined by (1.12), namely, $Z(t) := z_m(t)$, $t \in [t_m, t_{m+1})$, $m \in \mathbb{N}$, where for each $m \in \mathbb{N}$ the map $z_m(\cdot)$ is the solution of (2.55). Notice that for any initial state $x_0 \in \mathbb{R}^n$ of (2.1a), there exists $m_0 \in \mathbb{N}$ with $m_0 \geq 2$ such that

$$m_0 \geq |x_0| \tag{2.57}$$

Let $m \geq m_0$ and notice that, due to (2.55b), there exists a time $\bar{t} \in (t_{m-1}, t_{m+1}]$ such that $|z_m(t)| < \zeta_m$ for all $t \in [t_{m-1}, \bar{t})$. We claim that $\bar{t} = t_{m+1}$, namely

$$|z_m(t)| < \zeta_m, \ \forall t \in [t_{m-1}, t_{m+1}) \tag{2.58}$$

On the contrary, let $\bar{t} \in (t_{m-1}, t_{m+1})$ such that

$$|z_m(\bar{t})| = \zeta_m \text{ and } |z_m(t)| < \zeta_m, \ \forall t \in [t_{m-1}, \bar{t}) \tag{2.59}$$

therefore, by taking into account (2.55a),(2.56a) and (2.59), the map $z_m(\cdot)$ satisfies:

$$\dot{z}_m = F(t, z, y(t)) + \phi_m(t) P_m^{-1}(t) H'(t)(y(t) - H(t)z), \ \forall t \in [t_{m-1}, \bar{t}] \tag{2.60}$$

Define

$$e_m(t) := x(t) - z_m(t), \ \forall t \in [t_{m-1}, t_{m+1}] \tag{2.61}$$

Now, according to our Assumption A3, both A1 and A2 hold with $R := m$. Also, notice that, due to definition (2.43), the constant $\xi_m$ satisfies (2.30) with $R = m$. Since $|x_0| \leq m = R$, it follows by invoking (2.32) of Proposition 2.2 and (2.61):

$$|e_m(t)| < \beta(t_{m-1}, t_0, m)\sqrt{L} \exp\left\{\int_{t_{m-1}}^{t} -\bar{d}_m(s)ds\right\}, \ \forall t \in [t_{m-1}, \bar{t}] \tag{2.62}$$

By taking into account (2.43),(2.44) and (2.62) we get:

$$|e_m(t)| < \xi_m, \ \forall t \in [t_{m-1}, \bar{t}] \tag{2.63}$$

From (1.4),(2.56b),(2.61),(2.63) and by taking into account that $\bar{t} \leq t_{m+1}$ we deduce:

$$|z_m(\bar{t})| \leq |x(\bar{t})| + |e_m(\bar{t})| < \beta(\bar{t}, t_0, m) + \xi_m \leq \alpha(t_{m+1}, t_0, m) + \xi_m = \zeta_m \tag{2.64}$$

which contradicts (2.59), therefore (2.58) holds. The same arguments above assert that (2.62) and (2.63) hold for every $t \in [t_{m-1}, t_{m+1}]$. Notice next that (2.45) implies

$$\exp\left\{\int_{t_{m-1}}^{t} -\bar{d}_m(s)ds\right\} \leq \frac{1}{\beta(t_{m-1},t_0,m)m\sqrt{L}}, \quad \forall t \geq t_m \tag{2.65}$$

therefore, by taking into account (2.65) and the fact that the inequality in (2.62) holds for all $t \in [t_{m-1}, t_{m+1}]$, we have:

$$|e_m(t)| < \beta(t_{m-1},t_0,m)\sqrt{L}\exp\left\{\int_{t_{m-1}}^{t} -\bar{d}_m(s)ds\right\} \leq \frac{1}{m}, \quad \forall t \in [t_m, t_{m+1}] \tag{2.66}$$

Finally, we show that $\lim_{t \to \infty} |x(t) - Z(t)| = 0$; equivalently:

$$\forall \varepsilon > 0 \Rightarrow \exists T > 0 : |x(t) - Z(t)| < \varepsilon, \quad \forall t > T \tag{2.67}$$

Indeed, let $\varepsilon > 0$, $\bar{m} = \bar{m}(\varepsilon) \in \mathbb{N}$ with $\bar{m} \geq m_0$ such that $\frac{1}{\bar{m}} < \varepsilon$ and $T = T(\varepsilon) := t_{\bar{m}}$. Then by (2.47) it follows that for every $t > T$ there exists $m \geq \bar{m}$ such that $t_m \leq t < t_{m+1}$ and thus, since $m \geq m_0$ and by taking into account (1.12), (2.61) and (2.66), it follows:

$$|x(t) - Z(t)| = |x(t) - z_m(t)| = |e_m(t)| \leq \frac{1}{m} \leq \frac{1}{\bar{m}} < \varepsilon, \quad \forall t > T = T(\varepsilon)$$

The latter implies (2.67) and the proof of Proposition 2.3 is completed. ∎

## III. APPLICATION TO COMPOSITE SYSTEMS (1.2)

In this section we apply the main result of Section II (Proposition 2.3), to derive sufficient conditions for the solvability of the AC-SODP for a class of composite systems (1.2). We first need a preliminary technical result. Let $k, \ell, n_1, n_2 \in \mathbb{N}$ and consider the time-varying map:

$$J : \mathbb{R}_{\geq 0} \times \mathbb{R}^\ell \times \mathbb{R}^k \to \mathbb{R}^{(n_1+n_2) \times (n_1+n_2)}; \tag{3.1a}$$

$$J(t,q,y) := \begin{pmatrix} \mathsf{A}(q) & a(t,y)\mathsf{B}(t,y) \\ \mathsf{C}(q) & \mathsf{D}(t,q,y) \end{pmatrix} \tag{3.1b}$$

where $\mathsf{A} \in C^0(\mathbb{R}^\ell; \mathbb{R}^{n_1 \times n_1})$, $\mathsf{C} \in C^0(\mathbb{R}^\ell; \mathbb{R}^{n_2 \times n_1})$, $a \in C^1(\mathbb{R}_{\geq 0} \times \mathbb{R}^k; \mathbb{R})$, $\mathsf{B} \in C^1(\mathbb{R}_{\geq 0} \times \mathbb{R}^k; \mathbb{R}^{n_1 \times n_2})$ and $\mathsf{D} \in C^0(\mathbb{R}_{\geq 0} \times \mathbb{R}^\ell \times \mathbb{R}^k; \mathbb{R}^{n_2 \times n_2})$. Let $U$ and $W$ be nonempty subsets of $\mathbb{R}^n$ with $U \cap W \neq \emptyset$ and let $\Omega(\mathbb{R}_{\geq 0}, U)$ be a (nonempty) set of functions $y := y_{t_0, x_0} : [t_0, \infty) \to \mathbb{R}^k$ parameterized by $(t_0, x_0) \in \mathbb{R}_{\geq 0} \times U$. We make the following hypothesis for the pair $(\mathsf{B}, \mathsf{D})$ and the function $a(\cdot, \cdot)$ involved in (3.1):

**H2 •** *There exist constants $L > 1$ and $E > 0$ such that for every $t_0 \geq 0$ there exists a set-valued map $[t_0, \infty) \ni t \to Q(t) \subset \mathbb{R}^\ell$ satisfying the CP in such a way that for every $\bar{t}_0 \geq t_0$, $\tau_0 > 0$ and $y \in \Omega(\mathbb{R}_{\geq 0}, U) \cap C^0([t_0, \infty); \mathbb{R}^k)$, a time-varying symmetric matrix $P \in C^1([\bar{t}_0, \infty); \mathbb{R}^{n_2 \times n_2})$ and a function $d \in C^0([\bar{t}_0, \infty); \mathbb{R})$, both $\tau_0$-noncausal with respect to $\Omega(\mathbb{R}_{\geq 0}, U)$, can be found satisfying:*

$$P(t) \geq I_{n_2 \times n_2}, \quad \forall t \geq \bar{t}_0; \quad |P(\bar{t}_0)| \leq L; \tag{3.2a}$$

$$\int_{\bar{t}_0}^{t} d(s)ds > -E, \quad \forall t \geq \bar{t}_0; \quad \int_{\bar{t}_0}^{\infty} d(s)ds = \infty; \tag{3.2b}$$

$$e'P(t)\mathsf{D}(t,q,y(t))e + \tfrac{1}{2}e'\dot{P}(t)e \leq -d(t)e'P(t)e, \\ \forall t \geq \bar{t}_0, e \in \ker \mathsf{B}(t, y(t)), q \in Q(t), \text{ provided that (2.10) holds} \tag{3.2c}$$

- *We also assume that for every $t_0 \geq 0$ and $y \in \Omega(\mathbb{R}_{\geq 0}, U) \cap C^0([t_0, \infty); \mathbb{R}^k)$ it holds:*

$$a(t, y(t)) \neq 0, \text{ a.e. } t \geq t_0 \tag{3.3}$$

**Proposition 3.1.** *Consider the time-varying map $J(\cdot)$ as given by (3.1) and assume that H2 is fulfilled. Then for every $\varepsilon > 0$, $\bar{t}_0 \geq t_0 \geq 0$, $\tau > 0$ and $y \in \Omega(\mathbb{R}_{\geq 0}, U) \cap C^0([t_0, \infty); \mathbb{R}^k)$, there exist time-varying matrices $S \in C^1([\bar{t}_0, \infty); \mathbb{R}^{n_1 \times n_2})$, $T \in C^1([\bar{t}_0, \infty); \mathbb{R}^{n_1 \times n_1})$ and a function $\bar{d} \in C^0([\bar{t}_0, \infty); \mathbb{R})$, all being $\tau$-noncausal with respect to $\Omega(\mathbb{R}_{\geq 0}, U)$, such that, if we define*

$$\bar{P} := \begin{pmatrix} T & S \\ S' & P \end{pmatrix} \in C^1([\bar{t}_0, \infty); \mathbb{R}^{(n_1+n_2) \times (n_1+n_2)}) \tag{3.4}$$

*with $P(\cdot)$ as given in H2, the following properties hold:*

$$\bar{P}(t) \geq I_{(n_1+n_2) \times (n_1+n_2)}, \quad \forall t \geq \bar{t}_0; \quad |\bar{P}(\bar{t}_0)| \leq L; \tag{3.5a}$$

$$\int_{\bar{t}_0}^{t} \bar{d}(s) ds > -(E + \varepsilon), \quad \forall t \geq \bar{t}_0; \quad \int_{\bar{t}_0}^{\infty} \bar{d}(s) ds = \infty; \tag{3.5b}$$

$$(0, e') \bar{P}(t) J(t, q, y(t)) \begin{pmatrix} 0 \\ e \end{pmatrix} + \frac{1}{2}(0, e') \dot{\bar{P}}(t) \begin{pmatrix} 0 \\ e \end{pmatrix} \leq -\bar{d}(t)(0, e') \bar{P}(t) \begin{pmatrix} 0 \\ e \end{pmatrix}, \tag{3.5c}$$

$$\forall t \geq \bar{t}_0, e \in \mathbb{R}^{n_2}, q \in Q(t), \text{ provided that (2.10) holds}$$

***Proof.*** The proof is based on the result of Proposition 2.1. Notice that, due to (3.2c), the pair $(H, A)$ with $p = n_1$, $n = n_2$, $H(t, y) := \mathsf{B}(t, y)$ and $A(t, q, y) := \mathsf{D}(t, q, y)$ satisfies H1, with $P(\cdot)$, $d(\cdot)$ and $Q(\cdot)$ as given in (3.2). Let $\varepsilon$, $t_0$, $\bar{t}_0$, $\tau$ and $y(\cdot)$ as given in our statement, let $\tau_0 \in (0, \tau)$ and $\bar{\varepsilon} := \varepsilon / 2$. By exploiting (3.2b) and arguing as in proof of Corollary 2.1, a function $\hat{d} : [\bar{t}_0, \infty) \to \mathbb{R}$ can be found with

$$\hat{d}(t) < d(t), \quad \forall t \geq \bar{t}_0; \tag{3.6a}$$

$$\int_{\bar{t}_0}^{t} \hat{d}(s) ds > -(E + \bar{\varepsilon}) = -(E + \tfrac{\varepsilon}{2}), \quad \forall t \geq \bar{t}_0; \quad \int_{\bar{t}_0}^{\infty} \hat{d}(s) ds = \infty \tag{3.6b}$$

Since property H1 holds for $(A, H)$, it follows from (3.6a) and Proposition 2.1 that for any $\bar{\tau}_0 \in (\tau_0, \tau)$ there exists a function $\phi \in C^1([\bar{t}_0, \infty); \mathbb{R}_{>0})$, being $\bar{\tau}_0$-noncausal with respect to $\Omega(\mathbb{R}_{\geq 0}, U)$, such that

$$e' P(t) \mathsf{D}(t, q, y(t)) e + \tfrac{1}{2} e' \dot{P}(t) e \leq \phi(t) |\mathsf{B}(t, y(t)) e|^2 - \hat{d}(t) e' P(t) e, \tag{3.7}$$

$$\forall t \geq \bar{t}_0, e \in \mathbb{R}^{n_2}, q \in Q(t), \text{ provided that (2.10) holds}$$

Notice that, due to definitions (3.1) and (3.4), the desired (3.5c) is equivalent to

$$e' S'(t) a(t, y(t)) \mathsf{B}(t, y(t)) e + e' P(t) \mathsf{D}(t, q, y(t)) e + \tfrac{1}{2} e' \dot{P}(t) e \leq -\bar{d}(t) e' P(t) e, \tag{3.8}$$

$$\forall t \geq \bar{t}_0, e \in \mathbb{R}^{n_2}, q \in Q(t), \text{ provided that (2.10) holds}$$

It turns out, by taking into account (3.7) and (3.8), that, in order to show (3.5), it suffices to establish existence of mappings $T(\cdot)$, $S(\cdot)$ and $\bar{d}(\cdot)$ satisfying (3.5a),(3.5b), all being $\tau$-noncausal with respect to $\Omega(\mathbb{R}_{\geq 0}, U)$ in such a way that

$$e' S'(t) a(t, y(t)) \mathsf{B}(t, y(t)) e + \phi(t) |\mathsf{B}(t, y(t)) e|^2 - \hat{d}(t) e' P(t) e \leq -\bar{d}(t) e' P(t) e, \tag{3.9}$$

$$\forall t \geq \bar{t}_0, e \in \mathbb{R}^{n_2}, \text{ provided that (2.10) holds}$$

Let
$$S(t) := -\ell(t)a(t, y(t))\mathsf{B}(t, y(t)), \ t \geq \bar{t}_0 \tag{3.10}$$

for certain $\ell : [\bar{t}_0, \infty) \to \mathbb{R}_{\geq 0}$ satisfying

$$\ell(\bar{t}_0) = 0 \tag{3.11a}$$
$$\ell(t)a^2(t, y(t)) \leq \phi(t), \ \forall t \geq \bar{t}_0 \tag{3.11b}$$

yet to be determined and let

$$m(t) := \max\{|\mathsf{B}(t, y(t))w|^2 : w \in \mathbb{R}^{n_2}, |w| = 1\}, \ t \geq \bar{t}_0 \tag{3.12}$$

Obviously, $m(\cdot)$ is continuous and satisfies $m(t) \geq 0$ for every $t \geq \bar{t}_0$. We impose the following additional requirement for the desired $\bar{d}(\cdot)$:

$$\hat{d}(t) \geq \bar{d}(t), \ \forall t \geq \bar{t}_0 \tag{3.13}$$

By (3.9) and (3.10), it suffices to find $\ell(\cdot)$, $\bar{d}(\cdot)$ and $T(\cdot)$, all being $\tau$-noncausal with respect to $\Omega(\mathbb{R}_{\geq 0}, U)$, in such a way that (3.5a),(3.5b),(3.11) and (3.13) hold and further

$$(\phi(t) - \ell(t)a^2(t, y(t)))|\mathsf{B}(t, y(t))w|^2 \leq (\hat{d}(t) - \bar{d}(t))w'P(t)w$$
$$\forall t \geq \bar{t}_0, \ w \in \mathbb{R}^{n_2} : |w| = 1, \text{provided that (2.10) holds} \tag{3.14}$$

By taking into account the first inequality in (3.2a) and (3.12), it follows that, in order to show (3.14), it suffices to find $\ell(\cdot)$ and $\bar{d}(\cdot)$ satisfying (3.5b),(3.11) and (3.13), both being $\tau$-noncausal with respect to $\Omega(\mathbb{R}_{\geq 0}, U)$, in such a way that for $y \in \Omega(\mathbb{R}_{\geq 0}, U) \cap C^0([t_0, \infty); \mathbb{R}^k)$ as above it holds:

$$(\phi(t) - \ell(t)a^2(t, y(t)))m(t) \leq \hat{d}(t) - \bar{d}(t), \ \forall t \geq \bar{t}_0 \tag{3.15}$$

The desired (3.15) is a consequence of the following fact, whose proof constitutes a direct extension of the proof of Fact II in [49, page 1046]:

**Fact I:** Consider a nonempty set $\Omega(\mathbb{R}_{\geq 0}, U)$ of functions $y := y_{t_0, x_0} : [t_0, \infty) \to \mathbb{R}^k$ parameterized by $(t_0, x_0) \in \mathbb{R}_{\geq 0} \times U$, let $t_0 \geq 0$, $y \in \Omega(\mathbb{R}_{\geq 0}, U) \cap C^0([t_0, \infty); \mathbb{R}^k)$ and $\phi : [t_0, \infty) \to \mathbb{R}_{>0}$, $\zeta : [t_0, \infty) \to \mathbb{R}_{\geq 0}$, $\theta : [t_0, \infty) \to \mathbb{R}_{\geq 0}$ be continuous functions, whose values depend on $y(\cdot)$ and being $\tau_0$-noncausal with respect to $\Omega(\mathbb{R}_{\geq 0}, U)$. Moreover, assume that

$$\theta(t) \neq 0, \text{ a.e. } t \geq t_0 \tag{3.16}$$

Then for every $\varepsilon > 0$ and $\tau > \tau_0$, there exist a $C^1$ function $\ell : [t_0, \infty) \to \mathbb{R}_{\geq 0}$ and a $C^0$ function $h : [t_0, \infty) \to \mathbb{R}_{\geq 0}$, both $\tau$-noncausal with respect to $\Omega(\mathbb{R}_{\geq 0}, U)$, such that

$$(\phi(t) - \ell(t)\theta(t))\zeta(t) \leq h(t), \ \forall t \geq t_0 ; \tag{3.17a}$$
$$\phi(t) \geq \ell(t)\theta(t), \ \forall t \geq t_0 ; \tag{3.17b}$$
$$\int_{t_0}^{\infty} h(s)ds < \varepsilon ; \tag{3.17c}$$
$$\ell(t_0) = 0 \tag{3.17d}$$

By invoking requirements (3.11),(3.13) and (3.15), we may apply Fact I with $\phi(\cdot)$, $\zeta(\cdot) := m(\cdot)$ as given in (3.7) and (3.12), respectively, $\theta(\cdot) := a^2(\cdot, y(\cdot))$, (which due to (3.3) satisfies (3.16)), $t_0 := \bar{t}_0$, $\tau_0 := \bar{\tau}_0$ and $\varepsilon := \varepsilon/2$, in order to determine functions $\ell : [\bar{t}_0, \infty) \to \mathbb{R}_{\geq 0}$ and $h : [\bar{t}_0, \infty) \to \mathbb{R}_{\geq 0}$ being $\tau$-noncausal with respect to $\Omega(\mathbb{R}_{\geq 0}, U)$ and such that (3.17a)-(3.17d) hold. Then, if we define the $\tau$-noncausal with respect to $\Omega(\mathbb{R}_{\geq 0}, U)$ map $\bar{d} := \hat{d} - h$, inequality (3.13) is fulfilled, because $h(\cdot) \geq 0$. It also follows from (3.17a,b,d) that (3.15) as well as both (3.11a) and (3.11b) are satisfied. Moreover, (3.6b) and (3.17c) imply that $\bar{d}(\cdot)$ satisfies:

$$\int_{\bar{t}_0}^{t} \bar{d}(s)ds > -(E+\varepsilon), \ \forall t \geq \bar{t}_0 \ ; \ \int_{\bar{t}_0}^{\infty} \bar{d}(s)ds = \infty \tag{3.18}$$

thus (3.5b) holds as well. The proof is completed by the use of the following rather obvious fact:

**Fact II:** There exists a time-varying matrix $T \in C^1([\bar{t}_0, \infty); \mathbb{R}^{n_1 \times n_1})$, being $\tau$-noncausal with respect to $\Omega(\mathbb{R}_{\geq 0}, U)$, such that both inequalities in (3.5a) are satisfied with $\bar{P}(\cdot)$ as given by (3.4),(3.10) and (3.11a).

For reason of completeness an outline of proof of Fact II is given in the Appendix. ■

We are in a position to establish sufficient conditions for the solvability of the AC-SODP for the case of systems (1.2), where $f_1 : \mathbb{R}_{\geq 0} \times \mathbb{R}^{n_1} \to \mathbb{R}^{n_1}$ and $f_2 : \mathbb{R}_{\geq 0} \times \mathbb{R}^{n_1} \times \mathbb{R}^{n_2} \to \mathbb{R}^{n_2}$ are $C^0$ and locally Lipschitz on $x_1 \in \mathbb{R}^{n_1}$ and $(x_1, x_2) \in \mathbb{R}^{n_1} \times \mathbb{R}^{n_2}$, respectively. Also, assume that

$$G := a\mathsf{B} \tag{3.19}$$

for certain $a \in C^1(\mathbb{R}_{\geq 0} \times \mathbb{R}^{n_1}; \mathbb{R})$, $\mathsf{B} \in C^1(\mathbb{R}_{\geq 0} \times \mathbb{R}^{n_1}; \mathbb{R}^{n_1 \times n_2})$ and there exists a nonempty subset $M$ of $\mathbb{R}^{n_1} \times \mathbb{R}^{n_2}$ with $0 \in \text{cl}\,M$ such that (1.2a) is $M$-forward complete, namely, the solution $x(\cdot) := (x_1(\cdot), x_2(\cdot))$ of (1.2a) satisfies the estimation (1.4) for certain $\beta \in NNN$. Finally, define:

$$\Delta f_2(t, x_2, z_2; y) := f_2(t, y, x_2) - f_2(t, y, z_2); \tag{3.20}$$

$$H := (\underbrace{1,...,1}_{n_1}, \underbrace{0,...,0}_{n_2}) \tag{3.21}$$

and assume that there exists an integer $\ell \in \mathbb{N}$ and a map $\mathsf{D} \in C^0(\mathbb{R}_{\geq 0} \times \mathbb{R}^\ell \times \mathbb{R}^{n_1}; \mathbb{R}^{n_2 \times n_2})$ such that for every $R > 0$, $t_0 \geq 0$ and $\xi > 0$ there exists a set-valued map

$$[t_0, \infty) \ni t \to Q_R(t) := Q_{R, t_0, \xi}(t) \subset \mathbb{R}^\ell \tag{3.22}$$

satisfying the CP and in such a way that for every $t \geq t_0$ the following holds:

$$\begin{aligned}\exists q \in Q_R(t): \ \Delta f_2(t, x_2, z_2; y) = \mathsf{D}(t, q, y)(x_2 - z_2), \\ \text{provided that } y \in Y_R(t), \ x_2, z_2 \in \mathbb{R}^{n_2}, \ |x_2| \leq \beta(t, t_0, R) \text{ and } |x_2 - z_2| \leq \xi\end{aligned} \tag{3.23}$$

where $Y_R(t)$ is given in (2.4) with $H$ as defined by (3.21).

**Corollary 3.1.** *Let $M$ be a nonempty subset of $\mathbb{R}^{n_1} \times \mathbb{R}^{n_2}$ with $0 \in \text{cl}\,M$ and assume that (1.2a) is $M$-forward complete and its dynamics satisfy (3.19) and property (3.23). We further assume that there exists a constant $L > 1$ such that for every $R > 0$ there exists $\varepsilon_R > 0$ in such a way that for every $\bar{t}_0 \geq t_0 \geq 0$, $\tau_0 > 0$, $\xi > 0$ and $y \in Y(\mathbb{R}_{\geq 0}, M) \cap C^0([t_0, \infty); \mathbb{R}^{n_1})$, a time-varying symmetric matrix*

$$P_R := P_{R, t_0, \bar{t}_0, \tau_0, \xi, y} \in C^1([\bar{t}_0, \infty); \mathbb{R}^{n_2 \times n_2})$$

*and a function*

$$d_R := d_{R,t_0,\bar{t}_0,\tau_0,\xi,y} \in C^0([\bar{t}_0,\infty);\mathbb{R}),$$

*can be found, both $\tau_0$-noncausal with respect to $Y(\mathbb{R}_{\geq 0}, M)$, satisfying the following properties:*

$$P_R(t) \geq I, \ \forall t \geq \bar{t}_0; \ |P_R(\bar{t}_0)| \leq L; \tag{3.24a}$$

$$\int_{\bar{t}_0}^t d_R(s)ds > -\varepsilon_R, \ \forall t \geq \bar{t}_0; \ \int_{\bar{t}_0}^\infty d_R(s)ds = \infty; \tag{3.24b}$$

*and further*

$$e'P_R(t)\mathsf{D}(t,q,y(t))e + \tfrac{1}{2}e'\dot{P}_R(t)e \leq -d_R(t)e'P_R(t)e,$$
$$\forall t \geq \bar{t}_0, \ e \in \ker \mathsf{B}(t,y(t)), \ q \in Q_R(t), \text{provided that (2.7) holds} \tag{3.24c}$$

*where $Q_R(\cdot)$ and $\mathsf{D}(\cdot,\cdot,\cdot)$ are given in (3.22) and (3.23), respectively. Finally, for the function $a(\cdot)$ involved in (3.19) we assume that for every $t_0 \geq 0$ and $x_0 \in M$ it holds:*

$$a(t, x_1(t,t_0,x_0)) \neq 0 \text{ a.e. } t \geq t_0 \tag{3.25}$$

*Then the AC-SODP is solvable for (1.2) with respect to $Y(\mathbb{R}_{\geq 0}, M)$.*

***Proof.*** In order to establish our statement, it suffices to show that (1.2) satisfies A3, since according to Proposition 2.3 validity of A3 implies solvability of the AC-SODP. Define for each $t \geq 0$, $x = (x_1, x_2) \in \mathbb{R}^{n_1} \times \mathbb{R}^{n_2}$, $y \in \mathbb{R}^{n_1}$ and $q \in \mathbb{R}^\ell$:

$$F(t,x,y) := \begin{pmatrix} f_1(t,y) + a(t,y)\mathsf{B}(t,y)x_2 \\ f_2(t,y,x_2) \end{pmatrix}; \tag{3.26a}$$

$$A(t,q,y) := \begin{pmatrix} 0 & a(t,y)\mathsf{B}(t,y) \\ 0 & \mathsf{D}(t,q,y) \end{pmatrix} \tag{3.26b}$$

By taking into account (3.20),(3.23) and (3.26) it follows that for every $t \geq t_0$, $y \in Y_R(t)$ and $x, z \in \mathbb{R}^{n_1} \times \mathbb{R}^{n_2}$ with $|x| \leq \beta(t,t_0,R)$ and $|x-z| \leq \xi$ we have:

$$F(t,x,y) - F(t,z,y) = A(t,q,y)(x-z), \text{ for some } q \in Q_R(t) \tag{3.27}$$

thus for each $R > 0$ system (1.2) satisfies A1. We next show that for arbitrary $R > 0$ condition A2 is satisfied as well. Particularly, we show that for every $R > 0$ there exists $\bar{\varepsilon}_R > 0$ in such a way that for every $\bar{t}_0 \geq t_0 \geq 0$, $\tau > 0$, $\xi > 0$ and $y \in Y(\mathbb{R}_{\geq 0}, M) \cap C^0([t_0,\infty); \mathbb{R}^{n_1})$ there exist $\tau$-noncausal with respect to $Y(\mathbb{R}_{\geq 0}, M)$ mappings

$$\bar{P}_R \in C^1([\bar{t}_0,\infty); \mathbb{R}^{(n_1+n_2) \times (n_1+n_2)}), \ \bar{d}_R \in C^0([\bar{t}_0,\infty); \mathbb{R})$$

$\bar{P}_R(\cdot)$ being positive definite, in such a way that (2.6a,b) hold and further (2.6c) is fulfilled, provided that (2.7) holds with $L$, $A(\cdot,\cdot,\cdot)$ and $Q_R := Q_{R,t_0,\xi}$ as precisely given in (3.24a),(3.26b) and (3.22), respectively. We use the result of Proposition 3.1 to prove our statement. Let $R > 0$, $\xi > 0$, $\varepsilon_R$ as given in (3.24b) and define $\bar{\varepsilon}_R := 2\varepsilon_R$. Notice that $A(\cdot,\cdot,\cdot)$ as given by (3.26b) coincides with $J(\cdot,\cdot,\cdot)$ as given by (3.1) with $\mathsf{A}(q) := 0$, $\mathsf{C}(q) := 0$ and $a(\cdot,\cdot)$, $\mathsf{B}(\cdot,\cdot)$ and $\mathsf{D}(\cdot,\cdot,\cdot)$ as defined in (3.19) and (3.23). Then by invoking assumptions (3.24a,b,c) it follows that for every $\bar{t}_0 \geq t_0 \geq 0$ and $\tau_0 > 0$ conditions (3.2a,b,c) of H2 hold with $Q := Q_R(=Q_{R,t_0,\xi})$, $U := M$, $W := B_R$, $L$ as given in (3.24a), $E := \varepsilon_R$, $P := P_R$, $d := d_R$ and with $A(\cdot,\cdot,\cdot)$ as given by (3.26b). Also, (3.3) of H2 holds, due to (3.25). Therefore, all assumptions of Proposition 3.1 hold, hence, if we select $\varepsilon := \varepsilon_R$, then for every $\bar{t}_0 \geq t_0 \geq 0$, $\tau > 0$ and $y \in Y(\mathbb{R}_{\geq 0}, M) \cap C^0([t_0, \infty); \mathbb{R}^{n_1})$ there exist $\tau$-noncausal with respect to $Y(\mathbb{R}_{\geq 0}, M)$ time-varying mappings

$$S_R \in C^1([\bar{t}_0, \infty); \mathbb{R}^{n_1 \times n_2}), \quad T_R \in C^1([\bar{t}_0, \infty); \mathbb{R}^{n_1 \times n_1}), \quad \bar{d}_R \in C^0([\bar{t}_0, \infty); \mathbb{R})$$

such that, if we define:

$$\bar{P}_R := \begin{pmatrix} T_R & S_R \\ S_R' & P_R \end{pmatrix} \in C^1([\bar{t}_0, \infty); \mathbb{R}^{(n_1+n_2) \times (n_1+n_2)}) \tag{3.28}$$

with $P_R(\cdot)$ as given in (3.24a), the following properties are fulfilled:

$$\bar{P}_R(t) \geq I_{(n_1+n_2) \times (n_1+n_2)}, \quad \forall t \geq \bar{t}_0; \quad |\bar{P}_R(\bar{t}_0)| \leq L; \tag{3.29a}$$

$$\int_{\bar{t}_0}^t \bar{d}_R(s)ds > -(E+\varepsilon) = -2\varepsilon_R = -\bar{\varepsilon}_R, \quad \forall t \geq \bar{t}_0; \quad \int_{\bar{t}_0}^\infty \bar{d}_R(s)ds = \infty; \tag{3.29b}$$

and further

$$(0, e_2')\bar{P}_R(t)A(t,q,y)\begin{pmatrix} 0 \\ e_2 \end{pmatrix} + \frac{1}{2}(0, e_2')\dot{\bar{P}}_R(t)\begin{pmatrix} 0 \\ e_2 \end{pmatrix} \leq -\bar{d}_R(t)(0, e_2')\bar{P}_R(t)\begin{pmatrix} 0 \\ e_2 \end{pmatrix}, \tag{3.29c}$$

$\forall t \geq \bar{t}_0$, $e_2 \in \mathbb{R}^{n_2}$, $q \in Q_R(t)$, provided that (2.7) holds

By taking into account (3.21) and by setting $e := (e_1, e_2) \in \mathbb{R}^{n_1} \times \mathbb{R}^{n_2}$, (3.29c) becomes

$$e'\bar{P}_R(t)A(t,q,y(t))e + \tfrac{1}{2}e'\dot{\bar{P}}_R(t)e \leq -\bar{d}_R(t)e'\bar{P}_R(t)e, \tag{3.30}$$

$\forall t \geq \bar{t}_0$, $e \in \ker H$, $q \in Q_R(t)$, provided that (2.7) holds

By (3.29a), (3.29b) and (3.30) it follows that for arbitrary $R > 0$ all requirements of A2 hold as well. We conclude that (1.2) satisfies A3, therefore, according to Proposition 2.3, the AC-SODP is solvable for (1.2) with respect to $Y(\mathbb{R}_{\geq 0}, M)$. ∎

We illustrate the nature of Corollary 3.2 by the following examples.

**Example 3.1.** Consider the system

$$\begin{aligned}
\dot{x}_1 &= -x_1 + a(x_1)x_2 \\
\dot{x}_2 &= -x_1 a(x_1) + \frac{x_2}{x_1^2 + x_2^2 + x_3^2 + 1} \\
\dot{x}_3 &= -\frac{x_3}{x_1^2 + x_2^2 + 1}
\end{aligned} \tag{3.31a}$$

$$y = x_1 \tag{3.31b}$$

We assume that $a \in C^\omega(\mathbb{R}; \mathbb{R})$, ($C^\omega$ stands for analytic mappings), and there exists a nonempty $M_1$ subset of $\mathbb{R}$ with $0 \notin M_1$ and $0 \in \text{cl} M_1$ such that

$$\forall x \in M_1 \Rightarrow \exists m \in \mathbb{N}_0 : \frac{d^m}{dx^m}a(x) \neq 0 \tag{3.32}$$

Notice that system (3.31) has the form (1.2) with

$$f_1 := -x_1, \quad f_2 = \begin{pmatrix} f_2^1 \\ f_2^2 \end{pmatrix} := \begin{pmatrix} -x_1 a(x_1) + \dfrac{x_2}{x_1^2 + x_2^2 + x_3^2 + 1} \\ -\dfrac{x_3}{x_1^2 + x_2^2 + 1} \end{pmatrix},$$

$G(\cdot)$ as defined by (3.19) with

$$\mathsf{B} := (1 \quad 0) \tag{3.33}$$

and $a(\cdot)$ as given above. We claim that (3.31) satisfies all conditions of Corollary 3.2 with $M = M_1 \times \mathbb{R}^2$, therefore the AC-SODP is solvable for (3.31) with respect to $Y(\mathbb{R}_{\geq 0}, M_1 \times \mathbb{R}^{n_2})$. Indeed, by evaluating the time-derivative $\dot{V}$ of $V := x_1^2 + x_2^2 + x_3^2$ along the trajectories of (3.31a) we get $\dot{V} \leq 2$ for all $(x_1, x_2, x_3) \in \mathbb{R}^3$, therefore (3.31a) is forward complete; particularly, its solution $x(t) := x(t, x_0)$ initiated from $x_0$ at time $t = 0$ satisfies:

$$|x(t, x_0)| \leq \beta(t, |x_0|), \quad \beta(t, R) := 2\sqrt{t} + \sqrt{2}R, \quad t, R \geq 0 \tag{3.34}$$

Next, for each $R, \xi > 0$ consider functions $\underline{\sigma}_{R,\xi}^{i,j}, \bar{\sigma}_{R,\xi}^{i,j} \in C^0(\mathbb{R}_{\geq 0}; \mathbb{R})$, $1 \leq i, j \leq 2$ such that the following hold for $1 \leq i, j \leq 2$:

$$\bar{\sigma}_{R,\xi}^{i,j}(t) \geq \max\left\{\dfrac{\partial f_2^i(y, x_2, x_3)}{\partial x_{j+1}} : |(y, x_2, x_3)| \leq 2\beta(t, R) + \xi\right\}, \quad \forall t \geq 0; \tag{3.35a}$$

$$\underline{\sigma}_{R,\xi}^{i,j}(t) \leq \min\left\{\dfrac{\partial f_2^i(y, x_2, x_3)}{\partial x_{j+1}} : |(y, x_2, x_3)| \leq 2\beta(t, R) + \xi\right\}, \quad \forall t \geq 0 \tag{3.35b}$$

Notice that

$$\dfrac{\partial f_2^2}{\partial x_3}(y, x_2, x_3) = -\dfrac{1}{y^2 + x_2^2 + 1} \leq -\dfrac{1}{y^2 + x_2^2 + x_3^2 + 1} \leq -\dfrac{1}{2\beta(t, R)^2 + \xi + 1}, \tag{3.36}$$

$$\forall t \geq 0, \ x = (y, x_2, x_3) \in \mathbb{R}^3 : |x| \leq 2\beta(t, R) + \xi$$

hence, we may select in (3.35a):

$$\bar{\sigma}_{R,\xi}^{2,2}(t) := -\dfrac{1}{16t + 8R^2 + \xi + 1}, \quad t \geq 0 \tag{3.37}$$

thus, due to (3.34) and according to definition (3.37), we have:

$$\bar{\sigma}_{R,\xi}^{2,2}(t) \geq -\dfrac{1}{2\beta(t, R)^2 + \xi + 1}, \quad \forall t \geq 0 \tag{3.38}$$

Next, define:

$$Q_R(t) := \{q = (q_1, q_2, q_3, q_4) \in \mathbb{R}^4 : \underline{\sigma}_{R,\xi}^{i,j}(t) \leq |q_{i+2(j-1)}| \leq \bar{\sigma}_{R,\xi}^{i,j}(t), 1 \leq i, j \leq 2\}, \ t \geq 0; \tag{3.39a}$$

$$\mathsf{D}(q) := \begin{pmatrix} q_1 & q_2 \\ q_3 & q_4 \end{pmatrix}, \quad q := (q_1, q_2, q_3, q_4) \in \mathbb{R}^4 \tag{3.39b}$$

By taking into account (3.35a,b) and (3.39a,b), it follows that for every $t \geq 0$ there exists $q \in Q_R(t)$ such that (3.23) holds, where $Y_R(t)$ is given in (2.4) with $H := (1 \quad 0 \quad 0)$ and $\beta(\cdot)$ as given in (3.34). We establish that the remaining hypotheses of Corollary 3.2 are also fulfilled. Particularly, we show, that for every pair of constants $L > 1$ and $R > 0$ there exists $\varepsilon_R > 0$ such that for every $\bar{t}_0 \geq 0$, $\xi > 0$ and $y \in Y(\mathbb{R}_{\geq 0}, \mathbb{R}^3) \cap C^0(\mathbb{R}_{\geq 0}; \mathbb{R})$ there exist a time-varying symmetric matrix $P_R \in C^1([\bar{t}_0, \infty); \mathbb{R}^{2 \times 2})$ and a function $d_R \in C^0([\bar{t}_0,\infty);\mathbb{R})$, both being causal with respect to $Y(\mathbb{R}_{\geq 0},\mathbb{R}^3)$, satisfying (3.24a,b,c) and (3.25). Notice that, due to (3.33), the desired (3.24c) is equivalent to

$$(0,e)P_R(t)\mathsf{D}(q)\begin{pmatrix}0\\e\end{pmatrix}+\frac{1}{2}(0,e)\dot{P}_R(t)\begin{pmatrix}0\\e\end{pmatrix}\leq -d_R(t)(0,e)P_R(t)\begin{pmatrix}0\\e\end{pmatrix},\tag{3.40}$$

$$\forall t \geq \bar{t}_0,\ e \in \mathbb{R},\ q \in Q_R(t)$$

Define:

$$P_R(t):=\begin{pmatrix}p_1(t)&p(t)\\p(t)&p_2(t)\end{pmatrix}\tag{3.41}$$

for certain $p_1, p_2, p \in C^1([\bar{t}_0,\infty);\mathbb{R})$, yet to be determined. By taking into account (3.39a,b) and (3.41) the desired condition (3.40) is equivalent to

$$p(t)q_2 + p_2(t)q_4 + \tfrac{1}{2}\dot{p}_2(t) \leq -d_R(t)p_2(t),\tag{3.42}$$

$$\forall t \geq \bar{t}_0,\ (q_2,q_4)\in\mathbb{R}^2: \underline{\sigma}_{R,\xi}^{1,2}(t)\leq q_2 \leq \bar{\sigma}_{R,\xi}^{1,2}(t),\ \underline{\sigma}_{R,\xi}^{2,2}(t)\leq q_4 \leq \bar{\sigma}_{R,\xi}^{2,2}(t)$$

Define $p_1(t):=L$, $p_2(t):=L$, $p(t):=0$ and $d_R(t):=\dfrac{1}{16t+8R^2+\xi+1}$, $t\geq \bar{t}_0$ and notice that, due to (3.37) we have:

$$L(q_4 + d_R(t))\leq 0,\ \forall t \geq \bar{t}_0,\ q_4 \in \mathbb{R}: q_4 \leq \bar{\sigma}_{R,\xi}^{2,2}(t)\tag{3.43}$$

Then, the desired (3.42) is a consequence of (3.43) and it is obvious that (3.24a) and (3.24b) hold. Particularly, for every $\varepsilon_R > 0$ the first inequality in (3.24b) is fulfilled, because $d_R(\cdot)\geq 0$ and

$$\int_{\bar{t}_0}^{\infty}d_R(t)dt = \int_{\bar{t}_0}^{\infty}\frac{1}{16t+8R^2+\xi+1}dt = \infty$$

Finally, (3.25) is an immediate consequence of (3.32) and analyticity of the dynamics of (3.31a).

**Example 3.2.** Consider the autonomous system

$$\dot{x}_1 = f(x_1) + a(x_1)Bx_2$$
$$\dot{x}_2 = \bar{a}(x_1)\bar{B}x_2,\ (x_1,x_2)\in\mathbb{R}^{n_1}\times\mathbb{R}^{n_2}\tag{3.44a}$$
$$y = x_1\tag{3.44b}$$

where $n_1 < n_2$, $f \in C^{\omega}(\mathbb{R}^{n_1};\mathbb{R}^{n_1})$, $a \in C^{\omega}(\mathbb{R}^{n_1};\mathbb{R})$, $\bar{a}\in C^{\omega}(\mathbb{R}^{n_1};\mathbb{R})$ and we make the following assumptions:

- $(\bar{B},B)\in\mathbb{R}^{n_2\times n_2}\times\mathbb{R}^{n_1\times n_2}$ is a detectable pair of constant matrices, in the sense that the linear system $\dot{x}=\bar{B}x$ with output $y = Bx$ is detectable; also assume that there exists a positive definite matrix $S\in\mathbb{R}^{n_2\times n_2}$ such that

$$S\bar{B}+\bar{B}'S \leq 0\tag{3.45}$$

- there exists a nonempty $M_1$ subset of $\mathbb{R}^{n_1}$ with $0 \in \mathrm{cl}M_1$ such that

$$\forall x_1 \in M_1 \Rightarrow \exists m \in \mathbb{N}_0 : f^m a(x_1) \neq 0\tag{3.46}$$

where $f^k a := D(f^{k-1}a)f$, $k=1,2,\ldots;\ f^0 a := a$.

- we assume that

$$\bar{a}(x_1) > 0, \ \forall x_1 \in \mathbb{R}^{n_1} \quad (3.47a)$$

and there is a constant $C > 0$ such that

$$C \geq |a(x_1)|, \ \forall x_1 \in \mathbb{R}^{n_1} \quad (3.47b)$$

• finally, assume that there exists a positive definite function $V \in C^1(\mathbb{R}^{n_1}; \mathbb{R}_{\geq 0})$ and constants $C_i > 0$, $i = 1, 2, 3, 4$ such that

$$C_1 |x_1|^2 \leq V(x_1) \leq C_2 |x_1|^2, \ |DV(x_1)| \leq C_3 |x_1|, \ DV(x_1)f(x_1) \leq -C_4 |x_1|^2, \ \forall x_1 \in \mathbb{R}^{n_1} \quad (3.48)$$

System (3.44) has the form (1.2) with $G(\cdot)$ as defined by (3.19) with $a(\cdot)$ and $\mathsf{B}(\cdot) = B$ as given in (3.44a) and satisfies all conditions of Corollary 3.1 with $M = M_1 \times \mathbb{R}^{n_2}$, therefore the AC-SODP is solvable for (3.44) with respect to $Y(\mathbb{R}_{\geq 0}, M_1 \times \mathbb{R}^{n_2})$. Indeed, notice first that (3.45),(3.46) and (3.48) guarantee that (3.44) is forward complete, particularly, it can be easily verified that exists a function $\beta \in N$ such that the trajectory $x(\cdot) = (x_1(\cdot), x_2(\cdot))$ of (3.44) satisfies:

$$|x(t)| \leq \beta(|x(0)|), \ \forall t \geq 0, \ (x_1(0), x_2(0)) \in \mathbb{R}^{n_1} \times \mathbb{R}^{n_2} \quad (3.49)$$

The details of establishment of (3.49) are left to the reader. Also, by defining $\mathsf{D}(y) := \bar{a}(y)\bar{B}$ and taking into account (3.20) and (3.44a), it follows that equality in (3.23) holds for all $(x_2, z_2, y) \in \mathbb{R}^{n_2} \times \mathbb{R}^{n_2} \times \mathbb{R}^{n_1}$. Next, we take into account detectability of the pair $(\bar{B}, B)$ which according to [46] is equivalent to the existence of a constant positive definite matrix $P \in \mathbb{R}^{n_1 \times n_1}$ and a real constant $c > 0$ such that

$$P > I \ ; \ e'P\bar{B}e \leq -ce'Pe, \ \forall e \in \ker B \quad (3.50)$$

(Notice that $\ker B \neq \varnothing$, since $n_1 < n_2$). Now, for each $R > 0$ define $\rho_R := \min\{\bar{a}(x_1), |x_1| \leq \beta(R)\}$, (where $\beta(\cdot)$ is the gain function involved in (3.49)), which, due to (3.47a), is strictly positive, and let $P_R := P$ and $d_R := c\rho_R$. It follows from (3.49) and (3.50) that (3.24a,b,c) hold with $L := |P| > 1$. Finally, due to (3.46) and analyticity of dynamics, we can easily show, that condition (3.25) is fulfilled as well.

## IV. TRIANGULAR SYSTEMS

In this section we use the results of Sections II and III to establish sufficient conditions for the solvability of the AC-SODP for triangular systems (1.3). We assume that

$$f_i \in C^1(\mathbb{R}_{\geq 0} \times \mathbb{R}^i; \mathbb{R}), \ i = 1, 2, ..., n \text{ and } a_i \in C^1(\mathbb{R}_{\geq 0} \times \mathbb{R}; \mathbb{R}), \ i = 1, 2, ..., n-1$$

The following proposition generalizes Proposition 3.1 in [49].

**Proposition 4.1.** *Let $M$ be a nonempty subset of $\mathbb{R}^n$ with $0 \in \mathrm{cl}\,M$. Suppose that (1.3a) is M-forward complete, namely, the solution $x(\cdot) := x(\cdot, t_0, x_0)$ of (1.3a) satisfies the estimation (1.4) for certain $\beta \in NNN$ and also assume that for every $t_0 \geq 0$, $x_0 \in M$ and $i = 1, 2, ..., n-1$ it holds:*

$$a_i(t, x_1(t, t_0, x_0)) \neq 0 \text{ a.e. } t \geq t_0 \quad (4.1)$$

*Then the AC-SODP is solvable for (1.3) with respect to $Y(\mathbb{R}_{\geq 0}, M)$.*

*Proof.* To establish our statement, it suffices by virtue of the result of Proposition 2.3 to show that condition A3 is fulfilled for system (1.3), namely, we prove that there exist an integer $\ell \in \mathbb{N}$, a map $A(\cdot,\cdot,\cdot)$ and a constant $L > 1$, in such a way that for each $R > 0$ both A1 and A2 hold. Let $R > 0$, $t_0 \geq 0$ and $\xi > 0$ and define:

$$F(t, x, y) := \begin{pmatrix} f_1(t, y) + a_1(t, y) \\ f_2(t, y, x_2) + a_2(t, y) \\ \vdots \\ f_n(t, y, x_2, ..., x_n) \end{pmatrix}, (t, x, y) \in \mathbb{R}_{\geq 0} \times \mathbb{R}^n \times \mathbb{R} \quad (4.2)$$

Let $\sigma_{R, t_0, \xi} \in N$ such that

$$\sigma_{R, t_0, \xi}(t) \geq \sum_{i=2}^{n} \sum_{j=2}^{i} \max\left\{ \left| \frac{\partial f_i(t, y, x_2, ..., x_n)}{\partial x_j} \right| : |(y, x_2, ..., x_n)| \leq 2\beta(t, t_0, R) + \xi \right\}, \forall t \geq t_0 \quad (4.3)$$

and consider the set-valued map $[t_0, \infty) \ni t \to Q_R(t) := Q_{R, t_0, \xi}(t) \subset \mathbb{R}^\ell$, $\ell := \frac{n(n+1)}{2}$ defined as

$$Q_R(t) := \{q = (q_{1,1}; q_{2,1}, q_{2,2}; ...; q_{n,1}, q_{n,2}, ..., q_{n,n}) \in \mathbb{R}^\ell : |q| \leq \sigma_{R, t_0, \xi}(t)\} \quad (4.4)$$

that obviously satisfies the CP. Also, for each $t \geq t_0$ let $Y_R(t)$ as given in (2.4) with

$$H := (1, 0, ..., 0) \quad (4.5)$$

From (4.2)-(4.5) and use of the mean-value theorem it follows that for every $t \geq t_0$, $y \in Y_R(t)$ and $x, z \in \mathbb{R}^n$ with $|x| \leq \beta(t, t_0, R)$ and $|x - z| \leq \xi$ we have:

$$F(t, x, y) - F(t, z, y) = A(t, q, y)(x - z), \text{ for some } q \in Q_R(t) \text{ with } q_{i,1} = 0, i = 1, 2, ..., n; \quad (4.6a)$$

$$A(t, q, y) := \begin{pmatrix} q_{1,1} & a_1(t, y) & 0 & \cdots & 0 \\ q_{2,1} & q_{2,2} & a_2(t, y) & \ddots & \vdots \\ \vdots & \vdots & \vdots & \ddots & 0 \\ q_{n-1,1} & q_{n-1,2} & q_{n-1,3} & & a_{n-1}(t, y) \\ q_{n,1} & q_{n,2} & q_{n,3} & \cdots & q_{n,n} \end{pmatrix} \quad (4.6b)$$

thus A1 is satisfied. Next, by using the result of Proposition 3.1 of the previous section and adopting the same induction approach used in the proof of Proposition 3.1 in [49], we show that, under (4.1), A2 is satisfied as well. Particularly, we prove that there exists a constant $L > 1$ such that for every $R > 0$ there exists $\varepsilon_R > 0$ in such a way that for every $\overline{t_0} \geq t_0 \geq 0$, $\tau > 0$, $\xi > 0$ and $y \in Y(\mathbb{R}_{\geq 0}, M) \cap C^0([t_0, \infty); \mathbb{R})$, a time-varying symmetric matrix $P_R \in C^1([\overline{t_0}, \infty); \mathbb{R}^{n \times n})$ and a function $d_R \in C^0([\overline{t_0}, \infty); \mathbb{R})$ can be found, both $\tau_0$-noncausal with respect to $Y(\mathbb{R}_{\geq 0}, M)$, satisfying all conditions (2.6a,b), as well as (2.6c) provided that (2.7) holds, with $A(\cdot, \cdot, \cdot)$ and $Q_R := Q_{R, t_0, \xi}$ as given by (4.6b) and (4.4), respectively. In order to prove this claim, we proceed by induction as follows. Let $L > 1$, $R > 0$, $\varepsilon_R > 0$, $\xi > 0$ and for $k = 2, 3, ..., n$ define:

$$H_k := (1, 0, ..., 0) \in \mathbb{R}^k, \ e := (e_1, e_2, ..., e_k) \in \mathbb{R}^k; \quad (4.7a)$$

$$A_k(t,q,y) := \begin{pmatrix} q_{n-k+1,n-k+1} & a_{n-k+1}(t,y) & 0 & \cdots & 0 \\ q_{n-k+2,n-k+1} & & & & \\ \vdots & & \boxed{A_{k-1}(t,q,y)} & & \\ q_{n,n-k+1} & & & & \end{pmatrix} \quad (4.7b)$$

where

$$A_1(t,q,y) := q_{n,n} \quad (4.7c)$$

**Claim:** Let $k \in \mathbb{N}$ with $2 \leq k \leq n$. Then for $L$, $R$, $\varepsilon_R$, $\xi$ as above,

$$\bar{\varepsilon}_R := \tfrac{k}{n}\varepsilon_R \quad (4.8)$$

and for every $\bar{t}_0 \geq t_0 \geq 0$, $\tau > 0$ and $y \in Y(\mathbb{R}_{\geq 0}, M) \cap C^0([t_0, \infty); \mathbb{R})$, there exist a time-varying symmetric matrix $P_{R,k} \in C^1([\bar{t}_0, \infty); \mathbb{R}^{k \times k})$ and a mapping $d_{R,k} \in C^0([\bar{t}_0, \infty); \mathbb{R})$, both being $\tau$-noncausal with respect to $Y(\mathbb{R}_{\geq 0}, M)$, such that (2.6a,b) hold, and further (2.6c) is fulfilled, provided that (2.7) holds, with $d_R := d_{R,k}$, $H := H_k$, $P_R := P_{R,k}$, $A := A_k$ and $Q_R := Q_{R,t_0,\xi}$, namely, we have:

$$P_{R,k}(t) \geq I_{k \times k}, \ \forall t \geq \bar{t}_0; \ |P_{R,k}(\bar{t}_0)| \leq L; \quad (4.9a)$$

$$\int_{\bar{t}_0}^t d_{R,k}(s)ds > -\bar{\varepsilon}_R, \ \forall t \geq \bar{t}_0; \ \int_{\bar{t}_0}^\infty d_{R,k}(s)ds = \infty; \quad (4.9b)$$

$$e'P_{R,k}(t)A_k(t,q,y(t))e + \tfrac{1}{2}e'\dot{P}_{R,k}(t)e \leq -d_{R,k}(t)e'P_{R,k}(t)e$$
$$\forall t \geq \bar{t}_0, \ e \in \ker H_k, \ q \in Q_R(t), \text{provided that (2.7) holds} \quad (4.9c)$$

Notice that for $k := n$ the above claim guarantees that all requirements of A2 are fulfilled as well, therefore the desired statement of our proposition is a consequence of Proposition 2.3. We establish our claim by induction. Consider first the case $k := 2$, namely, let

$$H_2 := (1,0), \ e := (e_1, e_2) \in \mathbb{R}^2 \quad (4.10a)$$

$$A_2(t,q,y) := \begin{pmatrix} q_{n-1,n-1} & a_{n-1}(t,y) \\ q_{n,n-1} & q_{n,n} \end{pmatrix} \quad (4.10b)$$

Notice that $A_2(\cdot,\cdot,\cdot)$ as given in (4.10b) coincides with $J(\cdot,\cdot,\cdot)$ as given by (3.1) with $n_1 = n_2 = 1$ and

$$\mathsf{A}(q) := q_{n-1,n-1}, \mathsf{B}(t,y) := 1, a(t,y) := a_{n-1}(t,y),$$
$$\mathsf{C}(q) := q_{n,n-1} \text{ and } \mathsf{D}(t,q,y) := q_{n,n} \quad (4.11)$$

Consider the constants $L$, $R$, $\varepsilon_R$ and $\xi$ as above. Define:

$$P_{R,1}(t) := L, \ t \geq \bar{t}_0 \quad (4.12a)$$

and let $d_{R,1} \in C^0([\bar{t}_0, \infty); \mathbb{R})$ be a nonnegative function satisfying

$$\int_{\bar{t}_0}^\infty d_{R,1}(s)ds = \infty \quad (4.12b)$$

Obviously then, property H2 holds for the pair (B,D) as given in (4.12) with $Q := Q_{R,t_0,\xi}$ as defined by (4.4), $U := M$, $W := B_R$, $\Omega(\mathbb{R}_{\geq 0}, U) := Y(\mathbb{R}_{\geq 0}, M)$, $E := \tfrac{1}{n}\varepsilon_R$, and $P := P_{R,1}$, $d := d_{R,1}$ as given in (4.12). Indeed, by (4.12a) and (4.12b), it follows that (3.2a) and (3.2b) of property H2 are fulfilled. Also, due to (4.1), condition (3.3) of H2 is fulfilled and, due to (4.11), we have in our case that $\ker \mathsf{B}(t,y) = \{0\}$, therefore (3.2c) of H2 holds trivially. It follows by Proposition 3.1 that for every $\bar{t}_0 \geq t_0 \geq 0$, $\tau > 0$ and $y \in Y(\mathbb{R}_{\geq 0}, M) \cap C^0([t_0, \infty); \mathbb{R})$ there exist $\tau$-noncausal with respect to $Y(\mathbb{R}_{\geq 0}, M)$ time-varying mappings

$$S_{R,1} \in C^1([\bar{t}_0, \infty); \mathbb{R}), \ T_{R,1} \in C^1([\bar{t}_0, \infty); \mathbb{R}), \ d_{R,2} \in C^0([\bar{t}_0, \infty); \mathbb{R})$$

such that, if we define:

$$P_{R,2} := \begin{pmatrix} T_{R,1} & S_{R,1} \\ S'_{R,1} & P_{R,1} \end{pmatrix} \in C^1([\bar{t}_0, \infty); \mathbb{R}^{2 \times 2}) \tag{4.13}$$

with $P_{R,1}(\cdot)$ as given by (4.12a), the following properties hold:

$$P_{R,2}(t) \geq I_{2 \times 2}, \ \forall t \geq \bar{t}_0; \ |P_{R,2}(\bar{t}_0)| \leq L; \tag{4.14a}$$

$$\int_{\bar{t}_0}^{t} d_{R,2}(s) ds > -\tfrac{2}{n}\varepsilon_R, \ \forall t \geq \bar{t}_0; \ \int_{\bar{t}_0}^{\infty} d_{R,2}(s) ds = \infty; \tag{4.14b}$$

$$(0, e_2) P_{R,2}(t) A_2(t, q, y(t)) \begin{pmatrix} 0 \\ e_2 \end{pmatrix} + \tfrac{1}{2}(0, e_2) \dot{P}_{R,2}(t) \begin{pmatrix} 0 \\ e_2 \end{pmatrix} \leq -d_{R,2}(t)(0, e_2) P_{R,2}(t) \begin{pmatrix} 0 \\ e_2 \end{pmatrix}, \tag{4.14c}$$

$$\forall t \geq \bar{t}_0, \ e_2 \in \mathbb{R}, \ q \in Q_R(t), \text{provided that (2.7) holds}$$

By taking into account (4.10a), inequality (4.14c) becomes:

$$e'P_{R,2}(t)A_2(t,q,y(t))e + \tfrac{1}{2}e'\dot{P}_{R,2}(t)e \leq -d_{R,2}(t)e'P_{R,2}(t)e, \tag{4.15}$$
$$\forall t \geq \bar{t}_0, \ e \in \ker H_2, \ q \in Q_R(t), \text{provided that (2.7) holds}$$

From (4.7),(4.9),(4.14a,b) and (4.15) it follows that all properties of our claim hold for $k := 2$. In order to complete the proof, it suffices to show that, if the induction hypothesis of the claim holds for certain $k \in \{2,3,...,n-1\}$, then it is fulfilled for $k := k+1$ as well. We therefore show that for $L$, $R$, $\varepsilon_R$, $\xi$ as above and $\bar{\varepsilon}_R := \tfrac{k+1}{n}\varepsilon_R$ and for every $\bar{t}_0 \geq t_0 \geq 0$, $\tau > 0$ and $y \in Y(\mathbb{R}_{\geq 0}, M) \cap C^0([t_0, \infty); \mathbb{R})$, there exist a time-varying symmetric matrix $P_{R,k+1} \in C^1([\bar{t}_0, \infty); \mathbb{R}^{(k+1) \times (k+1)})$ and a map $d_{R,k+1} \in C^0([\bar{t}_0, \infty); \mathbb{R})$, both being $\tau$-noncausal with respect to $Y(\mathbb{R}_{\geq 0}, M)$, such that (2.6a,b) hold and further (2.6c) is fulfilled, provided that (2.7) holds, with $H := H_{k+1}$, $A := A_{k+1}$ and $Q_R := Q_{R,t_0,\xi}$ as given by (4.7a),(4.7b) and (4.4), respectively and appropriate $d_R := d_{R,k+1}$, $P_R := P_{R,k+1}$, yet to be determined. Notice that the map

$$A_{k+1}(t,q,y) = \begin{pmatrix} q_{n-k,n-k} & a_{n-k}(t,y) & 0 & \cdots & 0 \\ q_{n-k+1,n-k} & & & & \\ \vdots & & \boxed{A_k(t,q,y)} & & \\ q_{n,n-k} & & & & \end{pmatrix} \tag{4.16}$$

coincides with $J(t,q,y)$ as given by (3.1) with $n_1 = 1$, $n_2 = k$ and

$$\mathsf{A}(q) := q_{n-k,n-k}, \ \mathsf{C}(q) := \begin{pmatrix} q_{n-k+1,n-k} \\ \vdots \\ q_{n,n-k} \end{pmatrix}, \ a(t,y) := a_{n-k}(t,y),$$
$$\mathsf{B}(t,y) := \underbrace{(1,0,...,0)}_{k} \text{ and } \mathsf{D}(t,q,y) := A_k(t,q,y) \tag{4.17}$$

Again, notice that property H2 holds for the pair $(\mathsf{B},\mathsf{D})$ as given in (4.17) with $Q := Q_{R,t_0,\xi}$ as defined by (4.4), $U := M$, $W := B_R$, $\Omega(\mathbb{R}_{\geq 0}, U) := Y(\mathbb{R}_{\geq 0}, M)$, $E := \bar{\varepsilon}_R = \tfrac{k}{n}\varepsilon_R$ and $P := P_{R,k}$, $d := d_{R,k}$ as given in (4.9). Indeed, by (4.9a,b,c) it follows that (3.2a,b,c) of property H2 are fulfilled. Also, (3.3) of H2 holds, due to (4.1). Therefore, all assumptions of Proposition 3.1 hold, hence, if we select $\varepsilon := \frac{1}{n}\varepsilon_R$, then for every $\bar{t}_0 \geq t_0 \geq 0$, $\tau > 0$ and $y \in Y(\mathbb{R}_{\geq 0}, M) \cap C^0([t_0, \infty); \mathbb{R})$ we can find $\tau$-noncausal with respect to $Y(\mathbb{R}_{\geq 0}, M)$ time-varying mappings

$$S_{R,k+1} \in C^1([\bar{t}_0, \infty); \mathbb{R}^k), \ T_{R,k+1} \in C^1([\bar{t}_0, \infty); \mathbb{R}), \ d_{R,k+1} \in C^0([\bar{t}_0, \infty); \mathbb{R})$$

such that, if we define:

$$P_{R,k+1} := \begin{pmatrix} T_{R,k} & S_{R,k} \\ S'_{R,k} & P_{R,k} \end{pmatrix} \in C^1([\bar{t}_0, \infty); \mathbb{R}^{(k+1)\times(k+1)}) \tag{4.18}$$

with $P_{R,k}(\cdot)$ as given in (4.9a), the following properties hold:

$$P_{R,k+1}(t) \geq I_{(k+1)\times(k+1)}, \ \forall t \geq \bar{t}_0; \ |P_{R,k+1}(\bar{t}_0)| \leq L; \tag{4.19a}$$

$$\int_{\bar{t}_0}^{t} d_{R,k+1}(s)ds > -(\bar{\varepsilon}_R + \varepsilon) = -\frac{k+1}{n}\varepsilon_R, \ \forall t \geq \bar{t}_0; \ \int_{\bar{t}_0}^{\infty} d_{R,k+1}(s)ds = \infty; \tag{4.19b}$$

$$(0, e'_2)P_{R,k+1}(t)A_{k+1}(t,q,y(t))\begin{pmatrix} 0 \\ e_2 \end{pmatrix} + \frac{1}{2}(0, e'_2)\dot{P}_{R,k+1}(t)\begin{pmatrix} 0 \\ e_2 \end{pmatrix} \leq -d_{R,k+1}(t)(0, e'_2)P_{R,k+1}(t)\begin{pmatrix} 0 \\ e_2 \end{pmatrix}, \tag{4.19c}$$

$$\forall t \geq \bar{t}_0, \ e_2 \in \mathbb{R}^k, \ q \in Q_R(t), \text{provided that (2.7) holds}$$

By taking into account (4.7a) inequality (4.19c) becomes:

$$e'P_{R,k+1}(t)A_{k+1}(t,q,y(t))e + \frac{1}{2}e'\dot{P}_{R,k+1}(t)e \leq -d_{R,k+1}(t)e'P_{R,k+1}(t)e,$$
$$\forall t \geq \bar{t}_0, \ e \in \ker H_{k+1}, \ q \in Q_R(t), \text{provided that (2.7) holds} \tag{4.20}$$

It follows from (4.19a),(4.19b) and (4.20) that the induction hypothesis holds for $k := k+1$ as well. Therefore, for arbitrary $R > 0$ all requirements of A1 and A2 are fulfilled. We conclude that (1.3) satisfies A3, thus, according to Proposition 2.3, the AC-SODP is solvable for (1.3) with respect to $Y(\mathbb{R}_{\geq 0}, M)$. ∎

**Example 4.1.** Consider the 3-dimensional system

$$\begin{aligned}
\dot{x}_1 &= x_2 \\
\dot{x}_2 &= a(t, x_1)x_3 + g_1(t, x_1, x_2) - \xi_1 x_2^3 \\
\dot{x}_3 &= g_2(t, x_1, x_2, x_3) - \xi_2 x_3^3, \ (x_1, x_2, x_3) \in \mathbb{R}^3
\end{aligned} \tag{4.21a}$$

$$y = x_1 \tag{4.21b}$$

where $\xi_1, \xi_2 \geq 0$, $a \in C^1(\mathbb{R}_{\geq 0} \times \mathbb{R}; \mathbb{R})$ and the mappings $g_1 : \mathbb{R}_{\geq 0} \times \mathbb{R}^2 \to \mathbb{R}$ and $g_2 : \mathbb{R}_{\geq 0} \times \mathbb{R}^3 \to \mathbb{R}$ are $C^0$ and locally Lipschitz on $(x_1, x_2) \in \mathbb{R}^2$ and $(x_1, x_2, x_3) \in \mathbb{R}^3$, respectively. We make the following additional assumptions on system (4.21):

$$a(t, x_1) \neq 0, \ \forall (t, x_1) \in \mathbb{R}_{\geq 0} \times (\mathbb{R} \setminus \{0\}) \tag{4.22a}$$

$$g_1(t, 0, 0) = 0, \ \forall t \geq 0 \tag{4.22b}$$

there exists a constant $C > 0$ in such a way that for every $t \geq 0$ it holds:

$$|g_1(t, x_1, x_2)| \leq C|(x_1, x_2)|, \ \forall (x_1, x_2) \in \mathbb{R}^2 \text{ away from zero} \tag{4.23a}$$

$$|g_1(t, x_1, x_2, x_3)| \leq C|(x_1, x_2, x_3)|, \ \forall (x_1, x_2, x_3) \in \mathbb{R}^3 \text{ away from zero} \tag{4.23b}$$

System (4.21) has the form (1.3) with $a_1(t, x_1) := 1$, $a_2(t, x_1) := a(t, x_1)$ and it can be easily verified that due to (4.22) and (4.23), all conditions of Proposition 4.1 are fulfilled with $M := \{(x_1, x_2, x_3) \in \mathbb{R}^3 : (x_1, x_2) \neq 0\}$, therefore the AC-SODP is solvable with respect to $Y(\mathbb{R}_{\geq 0}, M)$. For completeness we note that by checking the time-derivative $\dot{V}$ of $V(x_1, x_2, x_3) := x_1^2 + x_2^2 + x_3^2$ along the trajectories of (4.21a) and taking into account (4.23a,b) it follows that $\dot{V}(x) \leq \ell V(x)$ for $x \in \mathbb{R}^3$ away from zero, thus (4.21a) is forward complete. The rest details for the validity of (4.1) are left to the reader.

## V. SOLVABILITY OF THE STRONG ODP

In this section we derive sufficient conditions for the solvability of the S-(S)ODP. Statements and proofs of the main results of this section constitute modifications of the statements and proofs of the corresponding results in Sections II, III and IV. Consider again the system (2.1), under the same regularity assumptions for its dynamics $F(\cdot, \cdot, \cdot)$ and the output map $H(\cdot)$. We assume that there exists a closed nonempty subset $M$ of $\mathbb{R}^n$ with $0 \in M$ such that (2.1a) is $M$-forward complete, namely (1.4) holds for certain $\beta \in NNN$. Also, assume that there exists an integer $\ell \in \mathbb{N}$, a map $A(\cdot, \cdot, \cdot)$ as in (2.2) and constants $L > 1$ and $R > 0$ such that A1 holds and, instead of A2, the following stronger hypothesis is fulfilled:

**A2'** *Condition A2 holds with the additional requirement that for every $\bar{t}_0 \geq t_0 \geq 0$, $\xi > 0$ and $y \in Y(\mathbb{R}_{\geq 0}, M) \cap C^0([t_0, \infty); \mathbb{R}^k)$ the mappings*

$$[\bar{t}_0, \infty) \times M \ni (t, x_0) \to P_{R, t_0, \bar{t}_0, \xi, y(\cdot, t_0, x_0)}(t), d_{R, t_0, \bar{t}_0, \xi, y(\cdot, t_0, x_0)}(t) \quad (5.1)$$
$$\text{are continuous and causal with respect to } Y(\mathbb{R}_{\geq 0}, M)$$

We next show that, under A1 and A2', the S-ODP is solvable for (2.1). As in Section II we first need to establish a preliminary result that constitutes a stronger version of the result of Proposition 2.1. Consider again the pair $(H, A)$ of continuous mappings in (2.8a,b), let $U$ and $W$ be nonempty subsets of $\mathbb{R}^n$, $U$ being closed and $W$ being compact, with $U \cap W \neq \emptyset$ and let $\Omega(\mathbb{R}_{\geq 0}, U)$ be a (nonempty) set of functions $y := y_{t_0, x_0} : [t_0, \infty) \to \mathbb{R}^k$ parameterized by $(t_0, x_0) \in \mathbb{R}_{\geq 0} \times \mathbb{R}^n$ in such a way that for every fixed $t_0 \in \mathbb{R}_{\geq 0}$, the map

$$[t_0, \infty) \times U \ni (t, x_0) \to y_{t_0, x_0}(t) \text{ is continuous} \quad (5.2)$$

Instead of H1, we make the following stronger hypothesis:

**H1'** *Condition H1 holds with the additional requirement that for every $\bar{t}_0 \geq t_0 \geq 0$ and $y := y_{t_0, x_0} \in \Omega(\mathbb{R}_{\geq 0}, U) \cap C^0([t_0, \infty); \mathbb{R}^k)$ the mappings*

$$[\bar{t}_0, \infty) \times U \ni (t, x_0) \to P_{t_0, \bar{t}_0, y_{t_0, x_0}}(t), d_{t_0, \bar{t}_0, y_{t_0, x_0}}(t) \text{ are continuous} \quad (5.3)$$
$$\text{and causal with respect to } \Omega(\mathbb{R}_{\geq 0}, U)$$

The following result (Proposition 5.1) constitutes a stronger version of Proposition 2.1:

**Proposition 5.1.** *Consider the pair $(H, A)$ of the continuous mappings in (2.8a,b) and assume that H1' is fulfilled. Then for every $\bar{t}_0 \geq t_0 \geq 0$ and for every $\bar{d} := \bar{d}_{t_0, \bar{t}_0, y_{t_0, x_0}} \in C^0([\bar{t}_0, \infty); \mathbb{R})$ satisfying (2.11) and in such a way that*

$$[\bar{t}_0, \infty) \times U \ni (t, x_0) \to \bar{d}_{t_0, \bar{t}_0, y_{t_0, x_0}}(t) \text{ is continuous} \quad (5.4)$$
$$\text{and causal with respect to } \Omega(\mathbb{R}_{\geq 0}, U)$$

*there exists a function* $\phi = \phi_{t_0,\bar{t}_0} \in C^1([\bar{t}_0,\infty);\mathbb{R}_{>0})$, *being independent of* $y(\cdot)$, *thus causal with respect to* $\Omega(\mathbb{R}_{\geq 0},U)$, *in such a way that (2.12) holds.*

**Proof(Outline).** Let $(t_0,x_0) \in \mathbb{R}_{\geq 0} \times (U \cap W)$, $y = y_{t_0,x_0} \in \Omega(\mathbb{R}_{\geq 0}, U \cap W) \cap C^0([t_0,\infty);\mathbb{R}^k)$ and define the mappings $D_y(t,q,e,x_0)$, $K(t,x_0)$, $K^c(t,x_0)$ and $\omega(t,x_0)$ as in (2.13a,b),(2.15) and (2.17), respectively. By compactness of $U \cap W$, continuity of $P(\cdot)$, $d(\cdot)$ with respect to $(t,x_0) \in [t_0,\infty) \times U \cap W$ and taking into account that (2.9),(2.10),(2.11),(2.13a) and the fact that $P(\cdot)$ is positive definite, it follows that $K(\cdot,\cdot)$ satisfies (2.14). Then, by arguing as in the proof of Proposition 2.1, we may show that for every $T > \bar{t}_0$ it holds

$$\inf\{\omega(t,x_0) : (t,x_0) \in [\bar{t}_0,T] \times U \cap W\} > 0 \tag{5.5}$$

By considering the map $\bar{\omega}(t,x_0)$ as in (2.19) it follows from (5.5) that for any $T > \bar{t}_0$ there exists a constant $M := M(T)$ such that

$$M \geq \sup\{\bar{\omega}(t,x_0) : (t,x_0) \in [\bar{t}_0,T] \times U \cap W\} \tag{5.6}$$

Define:

$$C(t,x_0) := \sup\{\bar{\omega}(t,x_0)(|P_{t_0,\bar{t}_0,y_{t_0,x_0}}(t)\| A(t,q,y(t))| + \tfrac{1}{2}|\dot{P}_{t_0,\bar{t}_0,y_{t_0,x_0}}(t)| + |\bar{d}_{t_0,\bar{t}_0,y_{t_0,x_0}}(t)\| P_{t_0,\bar{t}_0,y_{t_0,x_0}}(t)|) : q \in Q(t)\} \tag{5.7}$$

By taking into account (5.3),(5.4),(5.6),(5.7), compactness of $U \cap W$ and recalling the compactness property for the map $Q(\cdot)$, a causal function $\phi \in C^1([\bar{t}_0,\infty);\mathbb{R}_{>0})$, being independent of $y \in \Omega(\mathbb{R}_{\geq 0}, U \cap W) \cap C^0([t_0,\infty);\mathbb{R}^k)$ can be found, such that, (instead of (2.23)), it holds:

$$\phi(t) > \sup\{C(t,x_0), x_0 \in U \cap W\}, \ \forall t \geq \bar{t}_0 \tag{5.8}$$

By using similar arguments as in the first part of the proof of Proposition 2.1, we can establish that $\phi(\cdot)$ satisfies the desired (2.12). ∎

The results of Corollary 5.1 and Proposition 5.2 below constitute the "causal" analogues of Corollary 2.1 and Proposition 2.2, respectively. Their proofs are based on the result of Proposition 5.1 and are left to the reader.

**Corollary 5.1.** *Consider the pair* $(H,A)$ *as given in (2.8a,b) with* $H(t,y) := H(t)$ *and* $A(t,q,y)$ *as involved in (2.1b) and (2.2), respectively. Suppose that A2' is fulfilled and consider the constants* $R$, $\varepsilon_R$, $\xi$, $\bar{t}_0 \geq t_0 \geq 0$ *and the mappings*

$$Q_R := Q_{R,t_0,\xi}, \ y \in Y(\mathbb{R}_{\geq 0}, M) \cap C^0([t_0,\infty);\mathbb{R}^k), \ P_R(\cdot) := P_{R,t_0,\bar{t}_0,\xi,y}(\cdot), d_R(\cdot) := d_{R,t_0,\bar{t}_0,\xi,y}(\cdot)$$

*as determined in A2'. Then for every* $\bar{\varepsilon}_R > \varepsilon_R$, *there exists a function* $\bar{d}_R \in C^0([\bar{t}_0,\infty);\mathbb{R})$, *which is causal with respect to* $Y(\mathbb{R}_{\geq 0}, M)$ *and satisfies (2.27a,b), and a function* $\phi_R \in C^1([\bar{t}_0,\infty);\mathbb{R}_{>0})$ *being independent of* $y(\cdot)$, *thus causal with respect to* $Y(\mathbb{R}_{\geq 0}, M)$, *such that (2.27c) holds.*

**Proposition 5.2.** *Consider the system (2.1) and let* $M$ *be a nonempty closed subset of* $\mathbb{R}^n$ *with* $0 \in M$ *such that system (2.1a) is M-forward complete. For the initial state* $x_0 \in M$ *of (2.1a) assume that* $|x_0| \leq R$ *for some known constant* $R > 0$ *and assume that properties A1 and A2' hold for the constant* $R$ *above and for certain* $L > 1$. *Then the S-ODP is solvable for (2.1) with respect to* $Y(\mathbb{R}_{\geq 0}, B_R \cap M)$.

Next, we provide sufficient conditions for the solvability of the S-SODP. We strengthen A3 as follows:

**A3'** *Assume that there exist an integer $\ell \in \mathbb{N}$, a map $A(\cdot,\cdot,\cdot)$ as in (2.2) and a constant $L > 1$, in such a way that for every $R > 0$, hypotheses A1 and A2' hold.*

By exploiting the results of Corollary 5.1 and Proposition 5.2 and arguing as in proof of Proposition 2.3 we can establish the following causal version of Proposition 2.3:

**Proposition 5.3.** *In addition to hypothesis of M-forward completeness for (2.1a), assume that system (2.1) satisfies A3'. Then the S-SODP is solvable for (2.1) with respect to $Y(\mathbb{R}_{\geq 0}, M)$.*

Next, we derive sufficient conditions for the solvability of the S-SODP for composite systems of the form (1.2). Consider again the time-varying map $J(\cdot)$ as defined by (3.1), let $U$ and $W$ be nonempty subsets of $\mathbb{R}^n$, where $U$ is closed, $W$ compact and satisfy $U \cap W \neq \emptyset$ and consider a (nonempty) set $\Omega(\mathbb{R}_{\geq 0}, U)$ of functions $y := y_{t_0, x_0} : [t_0, \infty) \to \mathbb{R}^k$ parameterized by $(t_0, x_0) \in \mathbb{R}_{\geq 0} \times U$ in such a way that for every fixed $t_0 \in \mathbb{R}_{\geq 0}$ condition (5.2) holds. We strengthen hypothesis H2 for the pair $(\mathsf{B}, \mathsf{D})$ involved in (3.1) as follows:

**H2'** *Condition H2 holds with the additional requirements that for every $\overline{t}_0 \geq t_0 \geq 0$ and $y = y_{t_0, x_0} \in \Omega(\mathbb{R}_{\geq 0}, U) \cap C^0([t_0, \infty); \mathbb{R}^k)$ property (5.3) holds and further:*

$$a(t, y_{t_0, x_0}(t)) \neq 0, \quad \forall t \geq t_0, \ x_0 \in U \tag{5.9}$$

The following Proposition 5.4 and Corollary 5.2 are the causal versions of Proposition 3.1 and Corollary 3.1, respectively:

**Proposition 5.4.** *Consider the time varying map $J(\cdot)$ as given by (3.1) and assume that H2' is fulfilled. Then for every $\varepsilon > 0$, $\overline{t}_0 \geq t_0 \geq 0$ and $y_{t_0, x_0} \in \Omega(\mathbb{R}_{\geq 0}, U) \cap C^0([t_0, \infty); \mathbb{R}^k)$ there exist*

$$S := S_{t_0, \overline{t}_0, y_{t_0, x_0}} \in C^1([\overline{t}_0, \infty); \mathbb{R}^{n_1 \times n_2}), \ T := T_{t_0, \overline{t}_0, y_{t_0, x_0}} \in C^1([\overline{t}_0, \infty); \mathbb{R}^{n_1 \times n_1}), \ \overline{d} := \overline{d}_{t_0, \overline{t}_0, y_{t_0, x_0}} \in C^0([\overline{t}_0, \infty); \mathbb{R}),$$

*such that for each fixed $\varepsilon$, $t_0$ and $\overline{t}_0$ the mappings*

$$[\overline{t}_0, \infty) \times U \ni (t, x_0) \to S_{t_0, \overline{t}_0, y_{t_0, x_0}}(t), T_{t_0, \overline{t}_0, y_{t_0, x_0}}(t), \overline{d}_{t_0, \overline{t}_0, y_{t_0, x_0}}(t) \tag{5.10}$$
$$\text{are continuous and causal with respect to } \Omega(\mathbb{R}_{\geq 0}, U)$$

*and in such a way, that, if we define*

$$\overline{P} := \overline{P}_{t_0, \overline{t}_0, y_{t_0, x_0}} = \begin{pmatrix} T & S \\ S' & P \end{pmatrix} \in C^1([\overline{t}_0, \infty); \mathbb{R}^{(n_1 + n_2) \times (n_1 + n_2)}) \tag{5.11}$$

*with $P(\cdot)$ as given in H2', then (3.5a,b,c) hold.*

***Proof(Outline).*** The proof is based on the result of Proposition 5.2 and constitutes a simpler version of the proof of Proposition 3.1. Consider again the pair $(H, A)$ with $H(t, y) := \mathsf{B}(t, y)$ and $A(t, q, y) := \mathsf{D}(t, q, y)$ as given in (3.1) and notice that, due to (3.2c), property H1' holds for the pair $(H, A) = (\mathsf{B}, \mathsf{D})$, where for each $\overline{t}_0 \geq t_0 \geq 0$ and $y = y_{t_0, x_0} \in \Omega(\mathbb{R}_{\geq 0}, U) \cap C^0([t_0, \infty); \mathbb{R}^k)$, the corresponding mappings $P(\cdot)$ and $d(\cdot)$ are given in (3.2). Let $\varepsilon$, $t_0$, $\overline{t}_0$ and $y_{t_0, x_0}(\cdot)$ as given in our statement. Then, by recalling (3.2b) and (5.10) and by applying the same arguments with those in the proof of Corollary 2.1, we can find a causal function $\hat{d} := \hat{d}_{t_0, \bar{t}_0, y_{t_0, x_0}} \in C^0([\bar{t}_0, \infty); \mathbb{R})$ satisfying (3.6a,b) and in such a way that the map

$$[\bar{t}_0, \infty) \times U \ni (t, x_0) \to \hat{d}_{t_0, \bar{t}_0, y_{t_0, x_0}}(t) \text{ is continuous} \tag{5.12}$$

Thus, since H1' holds for the pair $(H, A) = (\mathsf{B}, \mathsf{D})$, it follows from (3.6a) and by invoking Proposition 5.1 that there exists a function $\phi \in C^1([\bar{t}_0, \infty); \mathbb{R}_{>0})$ being independent of $y = y_{t_0, x_0} \in \Omega(\mathbb{R}_{\geq 0}, U) \cap C^0([t_0, \infty); \mathbb{R}^k)$, thus causal with respect to $\Omega(\mathbb{R}_{\geq 0}, U)$, in such a way that (3.7) holds. As in the proof of Proposition 3.1, in order to prove (3.5a,b,c), it suffices to establish existence of mappings $T(\cdot)$, $S(\cdot)$ and $\bar{d}(\cdot)$ satisfying (5.10) and in such a way that (3.5a,b) and (3.9) hold. In our case the construction of the desired mappings is much simpler and is based on the additional assumption (5.9). Indeed, first, define:

$$\bar{m}(t) := \max\{|\mathsf{B}(t, y_{t_0, x_0}(t))w|^2 : w \in \mathbb{R}^{n_2}, |w| = 1, x_0 \in U \cap W\}, \ t \geq \bar{t}_0 \tag{5.13}$$

Notice that $\bar{m}(\cdot)$ is continuous and independent of $y(\cdot)$, thus causal with respect to $\Omega(\mathbb{R}_{\geq 0}, U)$. Furthermore, by taking into account definition (3.12) of the function $m(\cdot)$, it follows by (5.13) that

$$\bar{m}(t) \geq m(t), \ \forall t \geq \bar{t}_0, \text{ provided that (2.10) holds} \tag{5.14}$$

As in proof of Proposition 3.1, in order to get the desired (3.9), we find causal mappings $\ell(\cdot)$ and $\bar{d}(\cdot)$ satisfying (3.11) and (3.13) in such a way that (3.15) holds, provided that (2.10) is fulfilled. By taking into account (3.11b) and (5.14), it suffices to find $\ell(\cdot)$ and $\bar{d}(\cdot)$ satisfying:

$$(\phi(t) - \ell(t)a^2(t, y(t))\bar{m}(t) \leq \hat{d}(t) - \bar{d}(t), \forall t \geq \bar{t}_0, \text{provided that (2.10) holds} \tag{5.15}$$

Without any loss of generality we may assume that $\bar{m}(\cdot) \neq 0$, hence, a constant $\delta \in (0, \varepsilon)$ can be determined with $\int_{\bar{t}_0}^{T} \phi(s)\bar{m}(s)ds \geq \frac{\delta}{2}$ for certain $T > \bar{t}_0$. Define:

$$\Delta t := \min\left\{T > \bar{t}_0 : \int_{\bar{t}_0}^{T} \phi(s)\bar{m}(s)ds \geq \frac{\delta}{2}\right\} - \bar{t}_0 \tag{5.16}$$

Also, let $\eta \in C^1([\bar{t}_0, \infty); \mathbb{R})$ with $\eta(\bar{t}_0) := 0$, $\eta(t) \in [0, 1]$ for $t \in [\bar{t}_0, \bar{t}_0 + \Delta t]$ and $\eta(t) := 1$ for $t \geq \bar{t}_0 + \Delta t$. By taking into account (5.9), we define for the given $y_{t_0, x_0} \in \Omega(\mathbb{R}_{\geq 0}, U) \cap C^0([t_0, \infty); \mathbb{R}^k)$:

$$\ell(t) := \frac{\eta(t)\phi(t)}{a^2(t, y_{t_0, x_0}(t))}, \ t \geq \bar{t}_0 \tag{5.17}$$

Finally, let:

$$\bar{d}(t) := \hat{d}(t) - (\phi(t) - \ell(t)a^2(t, y(t)))\bar{m}(t), \ t \geq \bar{t}_0 \tag{5.18}$$

It follows from (5.17) and (5.18) that $\ell(\cdot)$ and $\bar{d}(\cdot)$ are causal, satisfying (3.11) and (3.13), respectively, and (5.15) holds, hence (3.9) is established. It remains to show that (3.5a,b) and (5.10) are fulfilled as well. Indeed, by taking into account (3.6b) and (5.16)-(5.18), it follows that $\bar{d}(\cdot)$ satisfies (3.5b). Furthermore, due to (5.12), (5.17) and (5.18), the function $\bar{d}(\cdot)$ satisfies (5.10) as well. Let $S(\cdot)$ as given by (3.10), which, due to (5.17), is causal and satisfies (5.10). Finally, as in the case of Proposition 3.1 (Fact II), we can construct a causal matrix $T \in C^1([\bar{t}_0, \infty); \mathbb{R}^{n_1 \times n_1})$ satisfying (5.10) as well as both inequalities in (3.5a) with $\bar{P}(\cdot)$ as given by (5.11). The details are left to the reader. ∎

**Corollary 5.2.** *Let $M$ be a nonempty closed subset of $\mathbb{R}^{n_1} \times \mathbb{R}^{n_2}$ with $0 \in M$ and assume that system (1.2a) satisfies all conditions imposed in Corollary 3.1 with the additional requirement that for each fixed $R$, $t_0$, $\bar{t}_0$ and $\xi$ the mappings*

$$[\bar{t}_0, \infty) \times M \ni (t, x_0) \to P_{R, t_0, \bar{t}_0, \xi, y(\cdot, t_0, x_0)}(t), \ d_{R, t_0, \bar{t}_0, \xi, y(\cdot, t_0, x_0)}(t) \tag{5.19}$$
$$\text{are continuous and causal with respect to } Y(\mathbb{R}_{\geq 0}, M)$$

*Also, assume that, instead of (3.25), for every $t_0 \geq 0$ and $x_0 \in M$ it holds:*

$$a(t, x_1(t, t_0, x_0)) \neq 0, \ \forall t \geq t_0 \tag{5.20}$$

*Then the S-SODP is solvable for (1.2) with respect to $Y(\mathbb{R}_{\geq 0}, M)$.*

**Proof.** The proof is based on the same arguments with those used in proof of Corollary 3.1, with the only exception that here we may invoke, due to the stronger assumptions (5.19) and (5.20), the results of Proposition 5.3 and Proposition 5.4, instead of Proposition 2.3 and Proposition 3.1, respectively. ∎

**Proposition 5.5.** *Let $M$ be a nonempty closed subset of $\mathbb{R}^n$ with $0 \in M$. Consider system (1.3) under the regularity assumptions of Proposition 4.1 for the mappings involved in the dynamics of (1.3a). Suppose that (1.3a) is $M$-forward complete, namely, the solution $x(\cdot) := x(\cdot, t_0, x_0)$ of (1.3a) satisfies the estimation (1.4) for certain $\beta \in NNN$ and also assume that for every $t_0 \geq 0$, $x_0 \in M$ and $i = 1, 2, \ldots, n-1$ it holds:*

$$a_i(t, x_1(t, t_0, x_0)) \neq 0, \ \forall t \geq t_0 \tag{5.21}$$

*Then the S-SODP is solvable for (1.3) with respect to $Y(\mathbb{R}_{\geq 0}, M)$.*

**Proof.** The proof is again based on the same arguments with those applied in proof of Proposition 4.1. The only exception here is that we invoke Proposition 5.3 and Proposition 5.4, instead of Proposition 2.3 and Proposition 3.1, respectively. ∎

The following examples illustrate the nature of Corollary 5.2 and Proposition 5.5.

**Example 5.1.** Consider the system (3.31) under the same regularity assumptions for the mapping $a(\cdot)$ and, instead of (3.32), we impose the stronger condition that $a(x) \neq 0$ for all $x \in \mathbb{R}$. Then it can be easily verified that (3.31) satisfies all conditions of Corollary 5.2 with $M = \mathbb{R}^3$, therefore the S-SODP is solvable for (3.31) with respect to $Y(\mathbb{R}_{\geq 0}, M)$.

**Example 5.2.** Consider the system (3.44) under the same regularity assumptions for the mappings $f(\cdot)$, $a(\cdot)$, $\bar{a}(\cdot)$ and the same detectability assumption concerning the pair $(\bar{B}, B)$. Also assume that (3.45),(3.47a,b) and (3.48) are fulfilled and, instead of (3.46), we impose the stronger condition $a(x_1) \neq 0$ for all $x_1 \in \mathbb{R}^{n_1}$. Then, as in the case of Example 3.2, it can be easily checked that (3.44) satisfies all conditions of Corollary 5.2 with $M = \mathbb{R}^{n_1} \times \mathbb{R}^{n_2}$, therefore the S-SODP is solvable for (3.44) with respect to $Y(\mathbb{R}_{\geq 0}, M)$.

**Example 5.3.** Finally, consider the system (4.21) under the same regularity assumptions for the mappings $a(\cdot)$, $g_1(\cdot)$ and $g_2(\cdot)$. Also assume that (4.22b) and (4.23a,b) are fulfilled and, instead of (4.22a), we assume that $a(t, x_1) \neq 0$ for all $(t, x_1) \in \mathbb{R}_{\geq 0} \times \mathbb{R}$. Then (4.21) satisfies all conditions of Proposition 5.5 with $M = \mathbb{R}^3$, therefore the S-SODP is solvable for (4.21) with respect to $Y(\mathbb{R}_{\geq 0}, M)$.

# APPENDIX

**Proof of Fact II.** We establish the existence of a time-varying matrix $T \in C^1([\bar{t}_0, \infty); \mathbb{R}^{n_1 \times n_1})$, being $\tau$-noncausal with respect to $\Omega(\mathbb{R}_{\geq 0}, U)$, such that both inequalities in (3.5a) are satisfied with $\bar{P}(\cdot)$ as given by (3.4),(3.10) and (3.11a). First, define:

$$I_k := (\underbrace{0}_{n_1-k} \ I_{k \times k}), \ S_k(t) := I_k S(t), \ k = 1, 2, \ldots, n_1, \ t \geq \bar{t}_0 \tag{A.1}$$

where $S(\cdot)$ is determined by (3.10) and (3.11a). Obviously, by (A.1),(3.10) and (3.11a), we have:

$$S_k(\bar{t}_0) = 0, \ k = 1, 2, \ldots, n_1 \tag{A.2}$$

Also, define:

$$T_k(t) := \begin{pmatrix} \tau_{k,k}(t) & 0 & \cdots & 0 \\ 0 & \ddots & 0 & \vdots \\ \vdots & 0 & \ddots & 0 \\ 0 & \cdots & 0 & \tau_{1,1}(t) \end{pmatrix}, \ k = 1, 2, \ldots, n_1, \ t \geq \bar{t}_0 \tag{A.3}$$

$$P_k(t) := \begin{pmatrix} T_k(t) & S_k(t) \\ S_k'(t) & P(t) \end{pmatrix}, \ k = 1, 2, \ldots, n_1, \ P_0(t) := P(t), \ t \geq \bar{t}_0 \tag{A.4}$$

for certain $\tau_{k,k} \in C^1([\bar{t}_0, \infty); \mathbb{R})$, $k = 1, 2, \ldots, n_1$, yet to be determined. We proceed by induction to construct appropriate $\tau_{k,k}(\cdot)$, $k = 1, 2, \ldots, n_1$ in such a way that each $P_k(\cdot)$ satisfies both requirements of (3.5a). Particularly, by taking into account (3.2a), (A.2)-(A.4) and by applying elementary induction procedure, we may show that for every $k \in \mathbb{N}$ with $1 \leq k \leq n_1$ and $i = 1, 2, \ldots, k$ it holds:

$$\det(P_i(t) - I_{(n_2+i) \times (n_2+i)}) := (\tau_{i,i}(t) - 1) \det(P_{i-1}(t) - I_{(n_2+i-1) \times (n_2+i-1)}) + K_i(t), \ t \geq \bar{t}_0; \tag{A.5a}$$

where

$$\tau_{i,i}(t) := L - \frac{K_i(t)}{\det(P_{i-1}(t) - I_{(n_2+i-1) \times (n_2+i-1)})}, \ t \geq \bar{t}_0; \tag{A.5b}$$

each $K_i \in C^1([\bar{t}_0, \infty); \mathbb{R})$, $i = 1, 2, \ldots, k$ above is independent of $\tau_{j,j}(\cdot)$, $i \leq j \leq k$ and simultaneously, the following hold:

$$\det(P_i(t) - I_{(n_2+i) \times (n_2+i)}) > 0, \ \forall t \geq \bar{t}_0; \tag{A.5c}$$

$$K_i(\bar{t}_0) = 0 \tag{A.5d}$$

Notice that, due to (3.2a) and (A.5c), each $P_k(\cdot)$, $k \in \mathbb{N}$, $1 \leq k \leq n_1$ satisfies the first inequality of (3.5a). Moreover, by (A.2)-(A.4) it follows that

$$P_k(\bar{t}_0) = \begin{pmatrix} \tau_{k,k}(\bar{t}_0) & \cdots & 0 & \\ \vdots & \ddots & \vdots & 0 \\ 0 & \cdots & \tau_{1,1}(\bar{t}_0) & \\ & 0 & & P(\bar{t}_0) \end{pmatrix}$$

thus, by taking into account the second inequality of (3.2a),(A.5b) and (A.5d), it follows that for every vector $x = (x_1, x_2, \ldots, x_{n_1+k}) \in \mathbb{R}^{n_1+k}$ with $|x| = 1$ it holds:

$$\mid \overline{P}(\overline{t_0})x \mid = \mid (\tau_{k,k}(\overline{t_0})x_1,...,\tau_{1,1}(\overline{t_0})x_k, P(\overline{t_0})(x_{k+1},...,x_{k+n_2})')\mid =$$

$$\left(\sum_{i=1}^{k}\mid \tau_{k+1-i,k+1-i}(\overline{t_0})\mid^2\mid x_i\mid^2 + \mid P(\overline{t_0})(x_{k+1},...,x_{k+n_2})'\mid^2\right)^{\frac{1}{2}} \leq \left(\sum_{i=1}^{k}L^2\mid x_i\mid^2 + \mid L(x_{k+1},...,x_{k+n_2})'\mid^2\right)^{\frac{1}{2}} = L\mid x\mid$$

therefore $P_k(\cdot)$ satisfies the second inequality of (3.5a) as well. We conclude that $\overline{P}(\cdot) := P_{n_1}(\cdot)$ satisfies both inequalities of (3.5a).